\documentclass[11pt]{article}

\usepackage[english]{babel}
\usepackage[T1]{fontenc}
\usepackage[utf8]{inputenc}
\usepackage{lmodern}
\usepackage{amsmath,amssymb,amsthm,mathtools,mathrsfs}
\usepackage{esint}
\usepackage{geometry}
\usepackage{microtype}
\usepackage{enumitem}
\usepackage{booktabs}
\usepackage{csquotes}
\usepackage{xurl}
\usepackage[hidelinks]{hyperref}
\usepackage[nameinlink,noabbrev]{cleveref}
\usepackage[backend=biber,style=numeric,sorting=nyt]{biblatex}
\AtBeginBibliography{\setlength{\emergencystretch}{2em}}

\newtheorem{theorem}{Theorem}[section]

\newtheorem{proposition}{Proposition}[section]

\AddToHook{env/proposition/before}{\stepcounter{theorem}}

\newtheorem{lemma}{Lemma}[section]

\AddToHook{env/lemma/before}{\stepcounter{theorem}}

\newtheorem{corollary}{Corollary}[section]

\AddToHook{env/corollary/before}{\stepcounter{theorem}}

\theoremstyle{definition}
\newtheorem{definition}{Definition}[section]

\AddToHook{env/definition/before}{\stepcounter{theorem}}

\AddToHook{env/example/before}{\stepcounter{theorem}}

\theoremstyle{remark}

\AddToHook{env/remark/before}{\stepcounter{theorem}}

\theoremstyle{plain}

\AddToHook{env/question/before}{\stepcounter{theorem}}

\crefname{theorem}{Theorem}{Theorems}
\Crefname{theorem}{Theorem}{Theorems}
\crefname{proposition}{Proposition}{Propositions}
\Crefname{proposition}{Proposition}{Propositions}
\crefname{lemma}{Lemma}{Lemmas}
\Crefname{lemma}{Lemma}{Lemmas}
\crefname{corollary}{Corollary}{Corollaries}
\Crefname{corollary}{Corollary}{Corollaries}
\crefname{definition}{Definition}{Definitions}
\Crefname{definition}{Definition}{Definitions}
\crefname{example}{Example}{Examples}
\Crefname{example}{Example}{Examples}
\crefname{remark}{Remark}{Remarks}
\Crefname{remark}{Remark}{Remarks}
\crefname{question}{Question}{Questions}
\Crefname{question}{Question}{Questions}
\crefname{section}{Section}{Sections}
\Crefname{section}{Section}{Sections}
\crefname{equation}{Equation}{Equations}
\Crefname{equation}{Equation}{Equations}

\title{Global Two-Sided Estimates for Time-Fractional Diffusion Kernels on Hyperbolic Space}
\author{Vicente Vergara\footnote{Department of Mathematics, Faculty of Physical and Mathematical Sciences, University of Concepci\'on, Concepci\'on, Chile. \texttt{vvergaraa@udec.cl}}}
\date{}

\hypersetup{
  pdftitle={Global Two-Sided Estimates for Time-Fractional Diffusion Kernels on Hyperbolic Space},
  pdfauthor={Vicente Vergara},
  pdfsubject={Global two-sided kernel estimates for time-fractional diffusion in hyperbolic geometry},
  pdfkeywords={time-fractional diffusion, hyperbolic space, two-sided kernel estimates, heat kernel, inverse stable subordination, saddle-point analysis}
}

\begin{document}
\maketitle

\begin{abstract}
We study time-fractional diffusion on real hyperbolic space $\mathbb H^n$, $n\ge2$, with Caputo order $0<\alpha<1$. We establish global two-sided pointwise estimates for the homogeneous fundamental kernel $Z_{\mathbb H}(t,r)$ and the inhomogeneous (Duhamel) kernel $Y_{\mathbb H}(t,r)$, uniformly over all $t>0$ and positive geodesic distances $r$. By combining inverse-stable subordination with global bounds for the hyperbolic heat kernel, we distinguish a short-time local regime, a large-time interior region, and an exterior region governed by a saddle-point problem with parameter $r/t$. Both kernels have the same exterior exponential rate, while their amplitudes are related by an explicit saddle-point factor. The estimates are uniformly comparable on the interfaces between regimes, recover the Euclidean time-fractional singularities at small scales, and identify the influence of hyperbolic geometry on large-scale decay. The kernel $Z_{\mathbb H}$ also describes the one-time marginal density of hyperbolic Brownian motion under an independent inverse stable clock. This representation yields a radial limit law, moment and quantile asymptotics, and a change in the mean square displacement from order $t^\alpha$ at small times to order $t^{2\alpha}$ at large times.
\end{abstract}

\medskip
\noindent\textbf{Keywords:} time-fractional diffusion; hyperbolic space; two-sided kernel estimates; heat kernel; inverse stable subordination; saddle-point analysis.

\noindent\textbf{2020 Mathematics Subject Classification:} Primary 35R11; Secondary 35K08, 58J35, 60K50.

\section{Introduction}
\label{sec:intro}

Let $\mathbb H^n$ be the real hyperbolic space of dimension $n\ge2$, and let
$0<\alpha<1$. We consider the time-fractional diffusion equation
\[
\partial_t^\alpha u-\Delta_{\mathbb H^n}u=f,
\]
where $\partial_t^\alpha$ denotes the Caputo derivative. The homogeneous problem
and the Duhamel term are governed, respectively, by the radial kernels
$Z_{\mathbb H}(t,r)$ and $Y_{\mathbb H}(t,r)$, with
$r=d_{\mathbb H^n}(x,y)$. Our aim is to obtain two-sided pointwise estimates,
uniform for all $t>0$ and $r>0$, and to make explicit the interaction between
temporal memory and the diffusion mechanisms induced by hyperbolic geometry.
In particular, the estimates must reflect both the locally Euclidean
structure of the hyperbolic heat kernel and its large-scale geometry. The
subordination formulas also give a probabilistic interpretation of
$Z_{\mathbb H}$ and a normalized representation for $Y_{\mathbb H}$; these
complementary consequences are discussed in \cref{sec:properties}.

Even before time fractionality is introduced, the heat equation on $\mathbb H^n$
exhibits dynamics that differ markedly from the Euclidean ones. Exponential
volume growth and the positive bottom of the spectrum are reflected in the heat
kernel and its pointwise decay; for the general theory of heat kernels and their
relation to geometry and spectrum, see Davies \cite{Davies1989}. In the hyperbolic
setting, Davies and Mandouvalos obtained precise global two-sided estimates
\cite{DaviesMandouvalos1988}, and Anker and Ostellari extend the framework to noncompact symmetric
spaces \cite{AnkerOstellari2003}. Vázquez also established a particularly
transparent dynamical consequence: at large times the radial mass of the
classical heat kernel concentrates around the line $r=(n-1)s$, with width of
order $s^{1/2}$, and a one-dimensional Gaussian profile emerges in the moving
frame \cite{Vazquez2022}. This ballistic propagation has no analogue in ordinary
Euclidean diffusion.

Replacing the first-order time derivative by a Caputo derivative introduces
memory and changes the structure of the fundamental solution. Caputo's
regularized derivative originated in models of dissipation with memory
\cite{Caputo1967}, and Metzler and Klafter discuss fractional diffusion in
anomalous transport and nonexponential relaxation \cite{MetzlerKlafter2000}.
Eidelman and Kochubei study the resulting fundamental kernels and emphasize
subdiffusive scaling in heterogeneous media \cite{EidelmanKochubei2004}.
The subordination approach also has a stochastic origin: Meerschaert and
Scheffler analyze scaling limits of continuous-time random walks with
heavy-tailed waiting times \cite{MeerschaertScheffler2004}, Baeumer and
Meerschaert establish stochastic representations of fractional Cauchy problems
\cite{BaeumerMeerschaert2001}, and Meerschaert, Nane, and Vellaisamy treat
fractional Cauchy problems on bounded domains
\cite{MeerschaertNaneVellaisamy2009}. Toaldo connects convolution-type time
operators with inverse subordinators and time-changed semigroups
\cite{Toaldo2015}; a broader account is given by Meerschaert and Sikorskii
\cite{MeerschaertSikorskii2019}.

From the analytic viewpoint, fundamental solutions of time-fractional diffusion
were studied early by Wyss and by Schneider--Wyss
\cite{Wyss1986,SchneiderWyss1989}, while Kochubei developed both the abstract
Cauchy problem and fractional-order diffusion
\cite{Kochubei1989,Kochubei1990}. The functional framework in Banach spaces and
its relation to semigroups was developed, among others, by Bazhlekova
\cite{Bazhlekova1998,Bazhlekova2001} and by Baeumer--Meerschaert through a
probabilistic subordination principle \cite{BaeumerMeerschaert2001}. In
$\mathbb R^d$, Kemppainen, Siljander, Vergara, and Zacher use fundamental
solutions to obtain optimal decay estimates and to identify critical dimensions
for time-fractional and other nonlocal-in-time diffusion equations
\cite{KemppainenSiljanderVergaraZacher2016}. For problems that are nonlocal both
in time and in space, Pozo--Vergara construct the associated kernels by
subordination and use them in representations of mild solutions and decay
estimates \cite{PozoVergara2019}. For elliptic operators with variable
coefficients, Eidelman--Kochubei construct and study the two kernels associated
with the initial datum and the inhomogeneous term
\cite{EidelmanKochubei2004}. In a direction especially close to ours,
Johnston--Kolokoltsov organize estimates for fractional Green functions from
known estimates for the kernel of the spatial operator, distinguishing which
information can be transferred globally or locally through subordination
\cite{JohnstonKolokoltsov2019}. In the Euclidean setting of stable
pseudodifferential generators, Solís--Vergara start from those estimates for
$Z$, derive the corresponding two-sided estimates for $Y$, and use both to
obtain $L^p$ properties and to analyze mild and global solutions of a
semilinear problem \cite{SolisVergara2022}. Mainardi, Mura, and Pagnini give a systematic account of the $M$--Wright
function and its asymptotic properties; this function describes the inverse
stable clock in our subordination formulas \cite{MainardiMuraPagnini2010}.

There is also prior work specifically on time-fractional diffusion in hyperbolic geometry.
Lao--Orsingher study planar hyperbolic Brownian motion and a Caputo time-fractional
version, obtaining an explicit representation of the law of the hyperbolic
distance; for order $1/2$ they also identify a composition with a reflected
Brownian time \cite{LaoOrsingher2007}. Garra--Maltese--Orsingher study, on
$\mathbb H^2$, a class of fractional derivatives with respect to another
function, obtain an explicit representation of the fundamental solution, and
give a time-change interpretation; choosing the identity time function recovers
the classical Caputo case \cite{GarraMalteseOrsingher2021}. These works provide explicit representations for fractional diffusion in
hyperbolic geometry. The present problem is to determine uniform, quantitative
pointwise bounds for both fundamental kernels, from the near-diagonal region to
the far field, including the interfaces between regimes. Starting from global
heat-kernel estimates on $\mathbb H^n$, we obtain two-sided bounds for
$Z_{\mathbb H}$ and $Y_{\mathbb H}$ in all dimensions $n\ge2$, with
explicitly identified local, interior, and exterior profiles.

The common starting point for both kernels is precisely this transfer by
subordination. We denote by $p_s(r)$ the hyperbolic heat kernel and by
$\ell_\alpha(t,s)$ the density of the inverse $\alpha$--stable clock; both will
be fixed precisely in \cref{sec:prelim}. Then
\[
Z_{\mathbb H}(t,r)
=
\int_0^\infty p_s(r)\ell_\alpha(t,s)\,\mathrm ds,
\qquad
Y_{\mathbb H}(t,r)
=
\frac{\alpha}{t}
\int_0^\infty s\,p_s(r)\ell_\alpha(t,s)\,\mathrm ds.
\]
Thus, after inserting the estimates for the hyperbolic heat kernel and the time
density, both kernels have the same exponential phase; the essential difference
is the additional operational-time factor $s$ in $Y_{\mathbb H}$. For
$Z_{\mathbb H}$, these are pointwise bounds for the marginal density of a
time-changed process, rather than transition bounds for a Markov semigroup.
Locally,
$\mathbb H^n$ resembles Euclidean space and the usual scales of time-fractional
diffusion reappear. Globally, however, the geometric drift and the radial escape
of the classical heat kernel compete with the distribution of the inverse
clock. Consequently, no single self-similar scale describes the whole $(t,r)$
plane uniformly.

\subsection{Main results}
\label{subsec:intro-main-results}

Our results are naturally organized into three regions: a short-time local
region, a large-time interior region, and an exterior region. This partition is
common to $Z_{\mathbb H}$ and $Y_{\mathbb H}$, but their profiles and amplitudes
are different. In each region we obtain quantitative two-sided pointwise
estimates, with profiles and rates completely determined by the parameters of
the problem. Throughout this summary we use the comparability convention fixed
in \cref{sec:prelim}.

For $Z_{\mathbb H}$, the local and interior regimes take the form
\[
Z_{\mathbb H}(t,r)
\asymp_{\alpha,n}
\begin{cases}
t^{-n\alpha/2}F_{\alpha,n}\!\left(r^2/t^\alpha\right),
&0<t\le1,\quad 0<r\le1,\\[2mm]
t^{-\alpha}G_{\mathbb H}(r),
&t>1,\quad 0<r\le t^\alpha,
\end{cases}
\]
where $F_{\alpha,n}$ and $G_{\mathbb H}$ are defined in
\eqref{eq:z-local-profile} and \eqref{eq:z-green-kernel}. In the exterior
region $r>\max\{1,t^\alpha\}$, we obtain the two-sided exterior estimate
\eqref{eq:z-exterior-estimate}. Its exponential rate is
\[
\exp\!\left[-r\left(\rho+J_\alpha(r/t)\right)\right],
\qquad \rho=\frac{n-1}{2},
\]
where the saddle point $y_\theta$ and the rate $J_\alpha(\theta)$ are defined
by \eqref{eq:z-saddle-equation} and \eqref{eq:z-rate-function}. Retaining
$J_\alpha$ rather than replacing it by its asymptotic expansions is essential
for preserving a uniform two-sided estimate throughout the exterior region.
The complete result is \cref{thm:z-global}.

For $Y_{\mathbb H}$, the local and interior regimes are
\[
Y_{\mathbb H}(t,r)
\asymp_{\alpha,n}
\begin{cases}
t^{\alpha-1-n\alpha/2}F^Y_{\alpha,n}\!\left(r^2/t^\alpha\right),
&0<t\le1,\quad 0<r\le1,\\[2mm]
t^{-\alpha-1}G_{\mathbb H}^{(2)}(r),
&t>1,\quad 0<r\le t^\alpha,
\end{cases}
\]
where $F^Y_{\alpha,n}$ and $G_{\mathbb H}^{(2)}$ are defined in
\eqref{eq:y-local-profile} and \eqref{eq:y-second-potential}. The additional
factor $s$ in the subordination formula shifts the local dimensional threshold
to $n=4$ and replaces the Green kernel by the potential kernel associated with $L^{-2}$. In the
exterior region, the full estimate is \eqref{eq:y-exterior-estimate}; the two
kernels share the same saddle point and exponential rate, and satisfy
\[
Y_{\mathbb H}(t,r)
\asymp_{\alpha,n}
\frac{r\,y_{r/t}}{t}\,Z_{\mathbb H}(t,r).
\]
The complete result is \cref{thm:y-global}.

The expressions obtained in the three regions are uniformly comparable on the interfaces
$t=1$, $r=1$, and $r=t^\alpha$. The local regime reproduces the singularities
of the Euclidean time-fractional kernels, whereas the large-time interior
region and the exterior rate record the hyperbolic geometry. In this sense,
the global estimates explicitly describe the transition between local Euclidean
behavior and large-scale geometric escape.

An additional consequence makes this transition visible at the probabilistic
level. For the classical heat equation on $\mathbb H^n$, mass propagates around
$r=(n-1)s$ \cite{Vazquez2022}. After the inverse-stable time change, if $R_t$
denotes the distance from the origin of the time-changed process and
$\mathcal{E}_1$ the inverse clock at time one, we prove, using
$\Longrightarrow$ to denote convergence in distribution, that
\[
\frac{R_t}{t^\alpha}
\Longrightarrow
(n-1)\mathcal{E}_1,
\qquad t\to\infty.
\]
Hence radial quantiles are of order $t^\alpha$, and the mean square displacement
changes scale from $t^\alpha$ as $t\downarrow0$ to $t^{2\alpha}$ as
$t\to\infty$. Unlike the Euclidean case, where the $t^\alpha$ law persists
globally, this crossover is a direct signature of hyperbolic ballistic
propagation. Mass normalizations, the probabilistic interpretation, and these
consequences are collected in \cref{sec:properties}.

\subsection{Organization of the paper}
\label{subsec:intro-organization}

In \cref{sec:prelim} we fix the geometric and fractional notation, recall the
global estimate for the hyperbolic heat kernel, and establish the subordinated
representations of $Z_{\mathbb H}$ and $Y_{\mathbb H}$. In
\cref{sec:z-estimates} we prove the estimates for $Z_{\mathbb H}$: the master
reduction and the local regime, the uniform analysis of the exterior saddle
point, and the large-time interior region, culminating in the global theorem.
\Cref{sec:y-estimates} reuses this structure for $Y_{\mathbb H}$,
identifies the $L^{-2}$ potential kernel, and transfers the exterior analysis
before closing its global estimate. Finally, \cref{sec:properties} collects the
geometric and probabilistic consequences: mass normalizations, time-change
representation, radial propagation, quantiles, and the mean-square-displacement
crossover.

\section{Preliminaries}
\label{sec:prelim}

Throughout the paper we fix $n\ge2$ and $0<\alpha<1$, and consider the real hyperbolic space $\mathbb H^n$ of constant sectional curvature $-1$. We denote by $d(\cdot,\cdot)$ the geodesic distance and by $\mathrm d\mu$ the Riemannian measure. For two-point kernels we use $r=d(x,y)$. When a pole $o\in\mathbb H^n$ is fixed and radial invariance is used, we employ the same symbol for
\[
r=d(x,o).
\]
The context will distinguish the two realizations of the same radial variable.

For two positive quantities $A$ and $B$, and a list of parameters $\mathcal P$, we write $A\asymp_{\mathcal P} B$
if there exist constants $0<c_{\mathcal P}\le C_{\mathcal P}<\infty$, depending only on $\mathcal P$, such that $c_{\mathcal P}B\le A\le C_{\mathcal P}B$.
Likewise, $A\lesssim_{\mathcal P}B$ means $A\le C_{\mathcal P}B$ and
$A\gtrsim_{\mathcal P}B$ means $A\ge c_{\mathcal P}B$. When the list
of parameters is clear from the context, we omit the subscript.

\subsection{Radial geometry and the heat kernel}

In geodesic polar coordinates around $o$, the volume element is
\[
\mathrm d\mu=(\sinh r)^{n-1}\,\mathrm dr\,\mathrm d\omega,
\]
where $\mathrm d\omega$ is the surface measure on $\mathbb S^{n-1}$. For a sufficiently regular radial function $v=v(r)$, the Laplace--Beltrami operator takes the form
\begin{equation}
\Delta_{\mathbb H^n}v
=
\partial_r^2v+(n-1)\coth r\,\partial_r v
=
(\sinh r)^{1-n}\partial_r\!\left((\sinh r)^{n-1}\partial_r v\right).
\label{eq:radial-laplacian}
\end{equation}

The radial form already anticipates the essential difference from the Euclidean case. Since
$\coth r\to1$ as $r\to\infty$, an outward radial drift with asymptotic speed
$n-1$ emerges far from the origin. Its large-scale probabilistic consequence
will be revisited in \cref{subsec:properties-radial}.

Let $p_s(x,y)$ be the kernel of the heat semigroup $e^{s\Delta_{\mathbb H^n}}$. By invariance under isometries, it depends only on the geodesic distance; we write
\[
p_s(x,y)=p_s(d(x,y)),
\qquad s>0.
\]
The normalization used throughout is
\[
\int_{\mathbb H^n}p_s(d(x,y))\,\mathrm d\mu(y)=1,
\qquad s>0.
\]
We also introduce the spectral constant
\[
\rho:=\frac{n-1}{2}.
\]
The number $\rho^2$ is the bottom of the spectrum of $-\Delta_{\mathbb H^n}$; this
quantity controls the spectral decay appearing in the heat kernel
\cite{Davies1989,Vazquez2022}.

The following result will be the basic geometric input for all subsequent subordination arguments.

\begin{proposition}[Davies--Mandouvalos]\label{prop:heat-kernel-dm}
For each $n\ge2$ there exist constants $0<c_n\le C_n<\infty$ such that, for every $s>0$ and $r\ge0$,
\begin{equation}
\begin{aligned}
c_n\,&s^{-n/2}(1+r)(1+r+s)^{(n-3)/2}
\exp\!\left(-\rho^2s-\rho r-\frac{r^2}{4s}\right)
\\
&\le p_s(r)\\
&\le C_n\,s^{-n/2}(1+r)(1+r+s)^{(n-3)/2}
\exp\!\left(-\rho^2s-\rho r-\frac{r^2}{4s}\right).
\end{aligned}
\label{eq:heat-kernel-dm}
\end{equation}
In particular,
\[
p_s(r)
\asymp_n
s^{-n/2}(1+r)(1+r+s)^{(n-3)/2}
\exp\!\left(-\rho^2s-\rho r-\frac{r^2}{4s}\right)
\]
uniformly in $(s,r)\in(0,\infty)\times[0,\infty)$.
\end{proposition}

The preceding estimate is a rewriting, in dimension $n$, of the global two-sided bound in \cite[Theorem~3.1]{DaviesMandouvalos1988}; for the broader setting of noncompact symmetric spaces see also \cite{AnkerOstellari2003}. We will repeatedly use the identity
\[
\rho^2s+\rho r+\frac{r^2}{4s}
=
\frac{(r+2\rho s)^2}{4s},
\]
which makes explicit the interaction between Gaussian diffusion and the hyperbolic geometric drift.

\subsection{Time-fractional calculus and subordination}
\label{subsec:fractional-subordination}

For $T>0$ and $f\in AC([0,T])$, where $AC([0,T])$ denotes the class of absolutely continuous functions, we use the regularized Caputo derivative of order $\alpha$, defined for almost every $t\in(0,T]$ by
\[
\partial_t^\alpha f(t)
=
\frac{1}{\Gamma(1-\alpha)}
\int_0^t
(t-\tau)^{-\alpha}f'(\tau)\,\mathrm d\tau.
\]
This class guarantees a finite initial trace $f(0)$ and makes the definition compatible with the classical Cauchy problem. This is the standard choice in time-fractional diffusion; its origin in dissipation models goes back to \cite{Caputo1967}, while early developments of fractional diffusion include \cite{Wyss1986,SchneiderWyss1989,Kochubei1990}.
Likewise, for $\beta>0$ we denote by
\[
E_{\alpha,\beta}(z)
=
\sum_{k=0}^{\infty}\frac{z^k}{\Gamma(\alpha k+\beta)}
\]
the two-parameter Mittag--Leffler function, and write $E_\alpha=E_{\alpha,1}$. We denote by $\mathcal L_t$ the Laplace transform with respect to the variable $t$; for vector-valued Laplace transforms and their relation to resolvents and abstract Cauchy problems, see \cite{ArendtBattyHieberNeubrander2011}.

\begin{definition}\label{def:m-wright}
The $M$--Wright function of index $\alpha$ is defined by
\[
M_\alpha(z)
=
\sum_{k=0}^{\infty}
\frac{(-z)^k}{k!\,\Gamma(1-\alpha-\alpha k)},
\qquad z\ge0.
\]
\end{definition}

The function $M_\alpha$ is continuous and positive on $[0,\infty)$, satisfies
\[
M_\alpha(0)=\frac{1}{\Gamma(1-\alpha)},
\]
and its moments are given by
\[
\int_0^\infty z^\beta M_\alpha(z)\,\mathrm dz
=
\frac{\Gamma(\beta+1)}{\Gamma(1+\alpha\beta)},
\qquad \beta>-1.
\]
Moreover, if
\[
q:=\frac{1}{1-\alpha},
\qquad
m:=q-1=\frac{\alpha}{1-\alpha},
\qquad
\kappa:=\frac{\alpha-\frac12}{1-\alpha},
\]
and
\[
b_\alpha:=(1-\alpha)\alpha^{\alpha/(1-\alpha)},
\]
then
\[
M_\alpha(z)
\sim
A_\alpha z^\kappa e^{-b_\alpha z^q},
\qquad z\to\infty,
\]
with
\[
A_\alpha
=
\frac{\alpha^\kappa}{\sqrt{2\pi(1-\alpha)}}.
\]
These properties are collected, for instance, in \cite{MainardiMuraPagnini2010}. By positivity and continuity on compact sets, the preceding asymptotics imply the global two-sided estimate
\begin{equation}
M_\alpha(z)
\asymp_\alpha
(1+z)^\kappa
\exp\!\left(-b_\alpha z^q\right),
\qquad z\ge0.
\label{eq:m-wright-global}
\end{equation}

Let $(D_s)_{s\ge0}$ be an $\alpha$--stable subordinator, normalized by
\[
\mathbb E\!\left[e^{-\lambda D_s}\right]
=
e^{-s\lambda^\alpha},
\qquad \lambda>0.
\]
Its inverse clock is defined by
\[
\mathcal{E}_t:=\inf\{s>0:D_s>t\},
\qquad t\ge0.
\]
For $t>0$, the random variable $\mathcal{E}_t$ has a density with respect to $s$, which under our
normalization is given by
\begin{equation}
\ell_\alpha(t,s)
=
t^{-\alpha}M_\alpha(st^{-\alpha}),
\qquad t>0,\ s>0.
\label{eq:inverse-density}
\end{equation}
In particular,
\[
\mathbb P(\mathcal{E}_t\in\mathrm ds)=\ell_\alpha(t,s)\,\mathrm ds,
\qquad
\int_0^\infty \ell_\alpha(t,s)\,\mathrm ds=1,
\]
and
\[
\int_0^\infty s^\beta\ell_\alpha(t,s)\,\mathrm ds
=
\frac{\Gamma(\beta+1)}{\Gamma(1+\alpha\beta)}
t^{\alpha\beta},
\qquad \beta>-1.
\]
Its Laplace transform in the time variable is
\[
\int_0^\infty e^{-\lambda t}\ell_\alpha(t,s)\,\mathrm dt
=
\lambda^{\alpha-1}e^{-s\lambda^\alpha},
\qquad \lambda>0.
\]

We will also use the density of the stable subordinator. With the order of
variables adopted here,
\begin{equation}
\eta_\alpha(t,s)
=
\frac{\alpha s}{t}\ell_\alpha(t,s)
=
\alpha s\,t^{-\alpha-1}M_\alpha(st^{-\alpha}),
\label{eq:stable-density}
\end{equation}
so that, for fixed $s>0$,
\[
\mathbb P(D_s\in\mathrm dt)=\eta_\alpha(t,s)\,\mathrm dt,
\qquad
\int_0^\infty e^{-\lambda t}\eta_\alpha(t,s)\,\mathrm dt
=
e^{-s\lambda^\alpha}.
\]

We recall the distinction between direct and inverse subordination. If $-A$
generates a bounded $C_0$--semigroup $T_s=e^{-sA}$ on a Banach space and its
clock is changed by the subordinator $D$, Bochner subordination produces, for
$\sigma>0$, the new semigroup
\[
T_\sigma^{(\alpha)}f
:=
\int_0^\infty T_s f\,\eta_\alpha(s,\sigma)\,\mathrm ds
=
e^{-\sigma A^\alpha}f.
\]
This construction goes back to Bochner and was developed for semigroups of
operators by Phillips; for a modern presentation in terms of Bernstein
functions, see \cite{Bochner1949,Phillips1952,SchillingSongVondracek2012}.
The Caputo time-fractional evolution studied below uses, instead, the inverse
clock $\mathcal{E}_t$ and the density $\ell_\alpha(t,s)$. The corresponding
representation is therefore called inverse subordination and is not, strictly
speaking, a Bochner-subordinated semigroup: the family of laws of
$\mathcal{E}_t$ does not form a convolution semigroup. From the viewpoint of
abstract Volterra equations, the same evolution can be formulated as
subordination of resolvent families in the sense of Prüss; in the Caputo case,
the density of the inverse stable clock gives an explicit realization of this
principle. See \cite{Pruss1993,PozoVergara2019} and, for the connection with
fractional Cauchy problems,
\cite{BaeumerMeerschaert2001,MainardiMuraPagnini2010,MeerschaertNaneVellaisamy2009,JohnstonKolokoltsov2019}.

We now set
\[
L:=-\Delta_{\mathbb H^n}\ge0.
\]
The two operator families entering the fractional equation are
\[
Z_{\mathbb H}(t)=E_\alpha(-t^\alpha L),
\qquad
Y_{\mathbb H}(t)=t^{\alpha-1}E_{\alpha,\alpha}(-t^\alpha L).
\]
This formulation belongs to the abstract theory of fractional Cauchy problems
and resolvent families associated with semigroup generators; see, among others,
\cite{Kochubei1989,Bazhlekova1998,Bazhlekova2001,BaeumerMeerschaert2001}.
Thus, for sufficiently regular data, the Duhamel representation for
\[
\partial_t^\alpha u+Lu=f,
\qquad
u(0)=u_0,
\]
is
\[
u(t)
=
Z_{\mathbb H}(t)u_0
+
\int_0^t
Y_{\mathbb H}(t-\tau)f(\tau)\,\mathrm d\tau.
\]

\begin{proposition}[Subordination representations]\label{prop:subordination-kernels}
For every $t>0$ and $r>0$,
\begin{equation}
Z_{\mathbb H}(t,r)
=
\int_0^\infty
p_s(r)\ell_\alpha(t,s)\,\mathrm ds,
\label{eq:subordination-z}
\end{equation}
and
\begin{equation}
\begin{aligned}
Y_{\mathbb H}(t,r)
&=
\int_0^\infty
p_s(r)\eta_\alpha(t,s)\,\mathrm ds
\\
&=
\frac{\alpha}{t}
\int_0^\infty
s\,p_s(r)\ell_\alpha(t,s)\,\mathrm ds.
\end{aligned}
\label{eq:subordination-y}
\end{equation}
In particular, both kernels are nonnegative and their spatial masses satisfy
\begin{equation}
\int_{\mathbb H^n}
Z_{\mathbb H}(t,d(o,x))\,\mathrm d\mu(x)=1,
\qquad
\int_{\mathbb H^n}
Y_{\mathbb H}(t,d(o,x))\,\mathrm d\mu(x)
=
\frac{t^{\alpha-1}}{\Gamma(\alpha)}.
\label{eq:kernel-masses}
\end{equation}
Moreover, for $r>0$,
\begin{equation}
Y_{\mathbb H}(\cdot,r)
=
D_t^{\,1-\alpha,\mathrm{RL}}
Z_{\mathbb H}(\cdot,r),
\label{eq:rl-relation}
\end{equation}
where
\[
D_t^{\,1-\alpha,\mathrm{RL}}g(t)
=
\frac{\mathrm d}{\mathrm dt}
\left[
\frac{1}{\Gamma(\alpha)}
\int_0^t
(t-\tau)^{\alpha-1}g(\tau)\,\mathrm d\tau
\right].
\]
\end{proposition}

\begin{proof}
Identity \eqref{eq:subordination-z} is the subordination formula associated with the inverse stable clock. For \eqref{eq:subordination-y}, using \eqref{eq:stable-density} and the Laplace transform in $t$,
\[
\mathcal L_t
\left[
\int_0^\infty
p_s(r)\eta_\alpha(t,s)\,\mathrm ds
\right](\lambda)
=
\int_0^\infty
p_s(r)e^{-s\lambda^\alpha}\,\mathrm ds.
\]
At the operator level, the right-hand side is the radial kernel of
\[
\int_0^\infty e^{-s(L+\lambda^\alpha)}\,\mathrm ds
=
(\lambda^\alpha+L)^{-1},
\]
which coincides with the Laplace transform of
$t^{\alpha-1}E_{\alpha,\alpha}(-t^\alpha L)$.
The second equality in \eqref{eq:subordination-y} follows directly from \eqref{eq:stable-density}.

The first identity in \eqref{eq:kernel-masses} follows from mass conservation for the heat semigroup and the fact that $\ell_\alpha(t,\cdot)$ is a probability density. For the second,
\[
\frac{\alpha}{t}
\int_0^\infty
s\,\ell_\alpha(t,s)\,\mathrm ds
=
\frac{\alpha}{t}
\frac{t^\alpha}{\Gamma(1+\alpha)}
=
\frac{t^{\alpha-1}}{\Gamma(\alpha)}.
\]
Finally,
\[
\mathcal L_t[Z_{\mathbb H}(\cdot,r)](\lambda)
=
\lambda^{\alpha-1}(\lambda^\alpha+L)^{-1}(r),
\]
whereas
\[
\mathcal L_t[Y_{\mathbb H}(\cdot,r)](\lambda)
=
(\lambda^\alpha+L)^{-1}(r).
\]
For $r>0$, $Z_{\mathbb H}(t,r)\to0$ as $t\downarrow0$, so the Laplace transform of the Riemann--Liouville derivative of order $1-\alpha$ multiplies by $\lambda^{1-\alpha}$. This proves \eqref{eq:rl-relation}. See also \cite{EidelmanKochubei2004} for the relation between the two kernels in the inhomogeneous fractional problem.
\end{proof}

Combining \eqref{eq:m-wright-global} with \eqref{eq:inverse-density} yields the form that will be used repeatedly:
\[
\ell_\alpha(t,s)
\asymp_\alpha
t^{-\alpha}
\left(1+\frac{s}{t^\alpha}\right)^\kappa
\exp\!\left[
-b_\alpha\left(\frac{s}{t^\alpha}\right)^q
\right].
\]
Together with the Davies--Mandouvalos estimate in \eqref{eq:heat-kernel-dm}, this identity reduces the estimates for $Z_{\mathbb H}$ and $Y_{\mathbb H}$ to one-dimensional integrals with the same phase and different amplitudes.

\section{Two-sided estimates for \texorpdfstring{$Z_{\mathbb H}$}{Z H}}
\label{sec:z-estimates}

In this section we study the homogeneous kernel $Z_{\mathbb H}(t,r)$. Every constant implicit in $\asymp_{\alpha,n}$ depends only on $(\alpha,n)$, unless explicitly stated otherwise.

\subsection{Master reduction}
\label{subsec:z-master-reduction}

The subordination representation \eqref{eq:subordination-z}, the global estimate for the $M$--Wright function \eqref{eq:m-wright-global}, and the two-sided heat-kernel bound in \cref{prop:heat-kernel-dm} reduce the problem to a positive one-dimensional integral.

\begin{proposition}[Master reduction]\label{prop:z-master-reduction}
For every $t>0$ and $r>0$,
\begin{equation}
\begin{aligned}
Z_{\mathbb H}(t,r)
\asymp_{\alpha,n}{}&
 t^{-\alpha}(1+r)e^{-\rho r}
\\
&\times
\int_0^\infty
s^{-n/2}(1+r+s)^{(n-3)/2}
\left(1+\frac{s}{t^\alpha}\right)^\kappa
\exp\!\left(-\Phi_{t,r}(s)\right)\,\mathrm ds,
\end{aligned}
\label{eq:z-master-reduction}
\end{equation}
where
\begin{equation}
\Phi_{t,r}(s)
=
\rho^2s+\frac{r^2}{4s}
+b_\alpha\left(\frac{s}{t^\alpha}\right)^q,
\qquad s>0.
\label{eq:z-master-phase}
\end{equation}
\end{proposition}

\begin{proof}
By \eqref{eq:subordination-z} and \eqref{eq:inverse-density},
\[
Z_{\mathbb H}(t,r)
=
t^{-\alpha}\int_0^\infty p_s(r)M_\alpha(st^{-\alpha})\,\mathrm ds.
\]
Applying \eqref{eq:heat-kernel-dm} and \eqref{eq:m-wright-global}, and grouping the three exponential terms that depend on $s$, yields \eqref{eq:z-master-reduction}--\eqref{eq:z-master-phase}. The constants are uniform in $(t,r)\in(0,\infty)^2$ and depend only on $(\alpha,n)$.
\end{proof}

The phase \eqref{eq:z-master-phase} contains simultaneously the hyperbolic spectral cost, the radial Gaussian cost, and the time-fractional subordination cost. The different regions of the global theorem will arise from different balances among these three terms.

\subsection{Profiles governing the regimes}
\label{subsec:z-governing-profiles}

At short times and on bounded spatial scales, the dimensionless variable
\[
\Omega:=\frac{r^2}{t^\alpha}
\]
appears. For $\Omega>0$, define
\begin{equation}
F_{\alpha,n}(\Omega)
:=
\int_0^\infty
\nu^{-n/2}(1+\nu)^\kappa
\exp\!\left(
-\frac{\Omega}{4\nu}-b_\alpha\nu^q
\right)\,\mathrm d\nu.
\label{eq:z-local-profile}
\end{equation}
This is the local profile obtained by scaling $s=t^\alpha\nu$ in \eqref{eq:z-master-reduction}.

At large times, inside the hyperbolic interior region, the natural object is the Green kernel of $-\Delta_{\mathbb H^n}$,
\begin{equation}
G_{\mathbb H}(r)
:=
\int_0^\infty p_s(r)\,\mathrm ds,
\qquad r>0.
\label{eq:z-green-kernel}
\end{equation}
The separation between these two profiles reflects two different mechanisms: at small scales, the balance between the Gaussian term and the $M$--Wright tail dominates, whereas in the large-time interior region the subordinating density is essentially constant over the relevant part of the Green integral.

\subsection{Short time in the local region}
\label{subsec:z-short-time-local}

Throughout this subsection we assume
\[
0<t\le1,
\qquad
0<r\le1,
\]
and retain the variable
\[
\Omega=\frac{r^2}{t^\alpha}.
\]
Since $r\le1$ and $t\le1$, if
\[
K:=t^{-\alpha}\ge1,
\]
then $0<\Omega\le K$.

Making the change of variables $s=t^\alpha\nu$ in \eqref{eq:z-master-reduction}, we obtain
\begin{equation}
Z_{\mathbb H}(t,r)
\asymp_{\alpha,n}
t^{-n\alpha/2}
\int_0^\infty
W_{t,r}(\nu)\,
f_\Omega(\nu)\,\mathrm d\nu,
\label{eq:z-short-rescaled}
\end{equation}
where
\[
f_\Omega(\nu)
:=
\nu^{-n/2}(1+\nu)^\kappa
\exp\!\left(
-\frac{\Omega}{4\nu}-b_\alpha\nu^q
\right)
\]
and
\[
W_{t,r}(\nu)
:=
(1+r)e^{-\rho r}
(1+r+t^\alpha\nu)^{(n-3)/2}
e^{-\rho^2t^\alpha\nu}.
\]
There exist constants $0<c_n\le C_n<\infty$ such that
\[
0<W_{t,r}(\nu)\le C_n,
\qquad \nu>0,
\]
and
\[
c_n\le W_{t,r}(\nu)\le C_n,
\qquad 0<\nu\le K.
\]
The delicate point is to justify that truncating the local profile at $\nu=K$ preserves comparability uniformly as $t\uparrow1$.

\begin{lemma}[Uniform truncation of the local profile]\label{lem:z-local-truncation}
Let
\[
I_K(\Omega):=\int_0^K f_\Omega(\nu)\,\mathrm d\nu.
\]
There exist constants $0<c_{\alpha,n}\le C_{\alpha,n}<\infty$ such that
\begin{equation}
c_{\alpha,n}F_{\alpha,n}(\Omega)
\le
I_K(\Omega)
\le
F_{\alpha,n}(\Omega)
\le
C_{\alpha,n}I_K(\Omega)
\label{eq:z-local-truncation}
\end{equation}
for every $K\ge1$ and $0<\Omega\le K$.
\end{lemma}

\begin{proof}
The inequality $I_K(\Omega)\le F_{\alpha,n}(\Omega)$ is immediate. Write
\[
T_K(\Omega)
:=
F_{\alpha,n}(\Omega)-I_K(\Omega)
=
\int_K^\infty f_\Omega(\nu)\,\mathrm d\nu.
\]
Since $\nu\ge K\ge1$, the factor $e^{-\Omega/(4\nu)}$ is at most one, and the algebraic factors have polynomial growth. Hence there exist exponents $B,B_1$ and positive constants, depending only on $(\alpha,n)$, such that
\[
T_K(\Omega)
\le
C_{\alpha,n}\int_K^\infty
\nu^B e^{-b_\alpha\nu^q}\,\mathrm d\nu
\le
C_{\alpha,n}K^{B_1}e^{-c_\alpha K^q}.
\]

Fix $K_0>1$. If $1\le K\le K_0$, integration over $\nu\in[1/2,1]$ gives
\[
I_K(\Omega)\ge c_{\alpha,n,K_0}>0,
\qquad 0<\Omega\le K,
\]
while $T_K(\Omega)$ remains uniformly bounded. Thus
$F_{\alpha,n}(\Omega)\lesssim_{\alpha,n,K_0} I_K(\Omega)$ in this range.

Now consider $K\ge K_0$. If $0<\Omega\le1$, the same interval $[1/2,1]$ provides a uniform lower bound for $I_K(\Omega)$, while the preceding bound for $T_K(\Omega)$ can be made arbitrarily small by increasing $K_0$.

It remains to consider $1\le\Omega\le K$. Let
\[
\nu_\Omega
:=
\left(
\frac{\Omega}{4\alpha^m}
\right)^{\frac{1-\alpha}{2-\alpha}}.
\]
Since
\[
\gamma:=\frac{1-\alpha}{2-\alpha}<1,
\qquad
\nu_\Omega\asymp_\alpha \Omega^\gamma,
\]
we may choose $K_0$ so that $2\nu_\Omega\le K$ for every
$K\ge K_0$ and $1\le\Omega\le K$. On
$\nu\in[\nu_\Omega,2\nu_\Omega]$, the amplitude is comparable to a power of $\Omega$, and
\[
\frac{\Omega}{4\nu}+b_\alpha\nu^q
\le
C_\alpha \Omega^{1/(2-\alpha)}.
\]
Consequently, there exist $A,C,c>0$, depending only on $(\alpha,n)$, such that
\[
I_K(\Omega)
\ge
c\,\Omega^{-A}
\exp\!\left(-C\Omega^{1/(2-\alpha)}\right)
\ge
c\,K^{-A}
\exp\!\left(-C K^{1/(2-\alpha)}\right).
\]
Since
\[
q=\frac{1}{1-\alpha}
>
\frac{1}{2-\alpha},
\]
the tail $T_K(\Omega)$ decays, as $K$ grows, faster than this lower bound. Increasing $K_0$ once more if necessary, we obtain
\[
T_K(\Omega)\le \frac12 I_K(\Omega).
\]
This completes the proof.
\end{proof}

\begin{proposition}[Short-time local regime]\label{prop:z-short-local}
For $0<t\le1$ and $0<r\le1$,
\begin{equation}
Z_{\mathbb H}(t,r)
\asymp_{\alpha,n}
t^{-n\alpha/2}
F_{\alpha,n}\!\left(\frac{r^2}{t^\alpha}\right).
\label{eq:z-short-local}
\end{equation}
\end{proposition}

\begin{proof}
The upper bound follows from \eqref{eq:z-short-rescaled} using
$W_{t,r}\le C_n$. For the lower bound, we restrict the integral to
$0<\nu\le K$, use $W_{t,r}\ge c_n$ on that interval, and apply
\cref{lem:z-local-truncation}.
\end{proof}

The profile $F_{\alpha,n}$ admits two explicit two-sided estimates, corresponding to
$\Omega\le1$ and $\Omega\ge1$.

\begin{proposition}[Two-sided asymptotics of the local profile]
\label{prop:z-local-profile-asymptotics}
For $0<\Omega\le1$,
\begin{equation}
F_{\alpha,n}(\Omega)
\asymp_{\alpha,n}
\begin{cases}
1+|\log\Omega|,
&
n=2,
\\[1mm]
\Omega^{1-n/2},
&
n\ge3.
\end{cases}
\label{eq:z-local-near}
\end{equation}
For $\Omega\ge1$,
\begin{equation}
F_{\alpha,n}(\Omega)
\asymp_{\alpha,n}
\Omega^{-\frac{n(1-\alpha)}{2(2-\alpha)}}
\exp\!\left(
-D_\alpha\Omega^{1/(2-\alpha)}
\right),
\label{eq:z-local-offdiag}
\end{equation}
where
\begin{equation}
D_\alpha
:=
(2-\alpha)
4^{-1/(2-\alpha)}
\alpha^{\alpha/(2-\alpha)}.
\label{eq:z-local-rate-constant}
\end{equation}
\end{proposition}

\begin{proof}
If $0<\Omega\le1$, the contribution of $\nu\ge1$ to
\eqref{eq:z-local-profile} is comparable to a constant. On $0<\nu\le1$,
the factors $(1+\nu)^\kappa$ and $e^{-b_\alpha\nu^q}$ are comparable to one, so
\[
F_{\alpha,n}(\Omega)
\asymp_{\alpha,n}
1+
\int_0^1
\nu^{-n/2}e^{-\Omega/(4\nu)}\,\mathrm d\nu.
\]
With the change of variables $u=\Omega/(4\nu)$,
\[
\int_0^1
\nu^{-n/2}e^{-\Omega/(4\nu)}\,\mathrm d\nu
=
4^{n/2-1}\Omega^{1-n/2}
\int_{\Omega/4}^{\infty}u^{n/2-2}e^{-u}\,\mathrm du.
\]
For $n\ge3$, the last integral remains comparable to a positive constant;
for $n=2$, it yields the logarithmic growth. This proves
\eqref{eq:z-local-near}.

Now suppose $\Omega\ge1$ and consider
\[
\phi_\Omega(\nu)
:=
\frac{\Omega}{4\nu}+b_\alpha\nu^q.
\]
Its unique minimizer is
\[
\nu_\Omega
=
\left(
\frac{\Omega}{4\alpha^m}
\right)^{\frac{1-\alpha}{2-\alpha}}.
\]
Define
\[
\lambda_\Omega
:=
\frac{\Omega}{4\nu_\Omega}
=
4^{-1/(2-\alpha)}
\alpha^{\alpha/(2-\alpha)}
\Omega^{1/(2-\alpha)}.
\]
Since $b_\alpha q=\alpha^m$, the change of variables $\nu=\nu_\Omega v$ gives
\[
\phi_\Omega(\nu_\Omega v)
=
\lambda_\Omega
\left(
\frac1v+\frac1qv^q
\right).
\]
The normalized phase has a unique nondegenerate minimum at $v=1$ and is coercive at both endpoints. Moreover, the normalized amplitude grows at most polynomially. The one-dimensional Laplace argument, uniformly for
$\lambda_\Omega\ge\lambda_1(\alpha)>0$, yields
\[
F_{\alpha,n}(\Omega)
\asymp_{\alpha,n}
\nu_\Omega^{1-n/2}(1+\nu_\Omega)^\kappa
\lambda_\Omega^{-1/2}
e^{-(1+1/q)\lambda_\Omega}.
\]
Since $\nu_\Omega\asymp_\alpha
\Omega^{(1-\alpha)/(2-\alpha)}$ and $1+1/q=2-\alpha$, simplifying the powers gives
\[
\nu_\Omega^{1-n/2}(1+\nu_\Omega)^\kappa
\lambda_\Omega^{-1/2}
\asymp_{\alpha,n}
\Omega^{-\frac{n(1-\alpha)}{2(2-\alpha)}},
\]
while
\[
(1+1/q)\lambda_\Omega
=
D_\alpha\Omega^{1/(2-\alpha)}.
\]
This proves \eqref{eq:z-local-offdiag}.
\end{proof}

The local powers in \eqref{eq:z-local-near} and the stretched-exponential rate in \eqref{eq:z-local-offdiag} agree with the structure of Euclidean time-fractional fundamental solutions; see \cite{Wyss1986,SchneiderWyss1989,Kochubei1990,EidelmanKochubei2004}. The global contribution here comes from the hyperbolic geometry outside this local scale.

Combining \cref{prop:z-short-local,prop:z-local-profile-asymptotics},
we obtain the explicit forms
\[
Z_{\mathbb H}(t,r)
\asymp_{\alpha,n}
\begin{cases}t^{-\alpha}
\left(
1+\left|\log\frac{r^2}{t^\alpha}\right|
\right),
&
n=2,\quad 0<r\le t^{\alpha/2},
\\[2mm]
t^{-\alpha}r^{2-n},
&
n\ge3,\quad 0<r\le t^{\alpha/2},
\end{cases}
\]
and
\[
\begin{aligned}
Z_{\mathbb H}(t,r)
\asymp_{\alpha,n}{}&
t^{-n\alpha/2}
\left(\frac{r^2}{t^\alpha}\right)^{
-\frac{n(1-\alpha)}{2(2-\alpha)}
}
\\
&\times
\exp\!\left[
-D_\alpha
\left(\frac{r^2}{t^\alpha}\right)^{1/(2-\alpha)}
\right],
\end{aligned}
\]
uniformly for
\[
0<t\le1,
\qquad
t^{\alpha/2}\le r\le1.
\]

\subsection{Parametrization of the exterior regime}
\label{subsec:z-exterior-saddle}

We now consider
\[
r\ge \max\{1,t^\alpha\}.
\]
Introduce
\[
\zeta:=\frac{r}{t^\alpha}\ge1,
\qquad
\theta:=\frac{r}{t}>0,
\qquad
s=ry.
\]
Since $r\ge1$, the algebraic factors in \eqref{eq:z-master-reduction} satisfy
\[
1+r\asymp r,
\qquad
1+r+ry\asymp r(1+y),
\]
uniformly for $y>0$. Consequently,
\begin{equation}
Z_{\mathbb H}(t,r)
\asymp_{\alpha,n}
 t^{-\alpha}r^{1/2}e^{-\rho r}
\int_0^\infty
A_\zeta(y)e^{-r h_\theta(y)}\,\mathrm dy,
\label{eq:z-exterior-reduction}
\end{equation}
where
\begin{equation}
A_\zeta(y)
:=
y^{-n/2}(1+y)^{(n-3)/2}(1+\zeta y)^\kappa
\label{eq:z-exterior-amplitude}
\end{equation}
and
\begin{equation}
h_\theta(y)
:=
\rho^2y+\frac{1}{4y}+b_\alpha\theta^m y^q.
\label{eq:z-exterior-phase}
\end{equation}

The phase $h_\theta$ is strictly convex on $(0,\infty)$, diverges at both endpoints, and has a unique minimizer, denoted by $y_\theta$. Since $b_\alpha q=\alpha^m$, the critical-point equation is
\begin{equation}
\frac{1}{4y_\theta^2}
=
\rho^2+\alpha^m\theta^m y_\theta^m.
\label{eq:z-saddle-equation}
\end{equation}
We define the exterior rate by
\begin{equation}
J_\alpha(\theta)
:=
h_\theta(y_\theta),
\qquad \theta>0.
\label{eq:z-rate-function}
\end{equation}
Thus the variable governing the transition of the critical point is $\theta=r/t$, whereas $\zeta=r/t^\alpha$ remains in the amplitude.

\begin{lemma}[Critical-point identities]\label{lem:z-saddle-identities}
For every $\theta>0$, the minimizer $y_\theta$ defined by \eqref{eq:z-saddle-equation} satisfies
\begin{equation}
J_\alpha(\theta)
=
\alpha\rho^2y_\theta
+
\frac{2-\alpha}{4y_\theta},
\label{eq:z-rate-identity}
\end{equation}
and, if
\[
\delta_\theta:=4\rho^2y_\theta^2\in(0,1],
\]
then
\begin{equation}
h_\theta''(y_\theta)
=
\frac{(q+1)-(q-1)\delta_\theta}{4y_\theta^3}.
\label{eq:z-curvature-identity}
\end{equation}
Moreover, with $y=y_\theta v$,
\begin{equation}
h_\theta(y_\theta v)
=
\frac{1}{4y_\theta}H_{\delta_\theta}(v),
\qquad
H_\delta(v)
:=
\delta v+\frac1v+\frac{1-\delta}{q}v^q,
\quad 0\le\delta\le1.
\label{eq:z-normalized-phase}
\end{equation}
In particular, $v=1$ is the unique minimizer of $H_\delta$ and
\[
2\le H_\delta''(1)\le q+1
\]
uniformly for $0\le\delta\le1$.
\end{lemma}

\begin{proof}
By \eqref{eq:z-saddle-equation},
\[
\alpha^m\theta^m y_\theta^m
=
\frac{1}{4y_\theta^2}-\rho^2.
\]
Since $b_\alpha=(1-\alpha)\alpha^m$, substituting this identity into $h_\theta(y_\theta)$ gives \eqref{eq:z-rate-identity}. Differentiating \eqref{eq:z-exterior-phase} twice and using the critical-point equation once more gives \eqref{eq:z-curvature-identity}.

Finally, multiplying $h_\theta(y_\theta v)$ by $4y_\theta$ and using
\[
4\rho^2y_\theta^2=\delta_\theta,
\qquad
4b_\alpha\theta^m y_\theta^{q+1}
=
\frac{1-\delta_\theta}{q},
\]
gives \eqref{eq:z-normalized-phase}. The second derivative is
\[
H_\delta''(v)
=
\frac{2}{v^3}
+(1-\delta)(q-1)v^{q-2}>0,
\]
and at $v=1$ it equals $2+(1-\delta)(q-1)$, which gives the stated uniform bound.
\end{proof}

The normalization \eqref{eq:z-normalized-phase} makes it possible to treat all values of the exterior parameter simultaneously. The next lemma contains the required uniformity.

\begin{lemma}[Uniform Laplace estimate in the exterior regime]\label{lem:z-uniform-laplace}
Let
\[
\Psi_\delta(v):=H_\delta(v)-H_\delta(1),
\qquad 0\le\delta\le1,
\]
and, for $\zeta\ge1$, $\theta>0$, write
\[
R_{\zeta,\theta}(v)
:=
\frac{A_\zeta(y_\theta v)}{A_\zeta(y_\theta)}.
\]
There exists a constant $C_{\alpha,n}\ge1$ such that, for every $\zeta\ge1$ and every $\theta>0$,
\begin{equation}
C_{\alpha,n}^{-1}\lambda^{-1/2}
\le
\int_0^\infty
R_{\zeta,\theta}(v)
e^{-\lambda\Psi_{\delta_\theta}(v)}\,\mathrm dv
\le
C_{\alpha,n}\lambda^{-1/2},
\qquad
\lambda\ge\frac{\rho}{2}.
\label{eq:z-uniform-laplace}
\end{equation}
Consequently, if $r\ge1$,
\begin{equation}
r^{1/2}
\int_0^\infty
A_\zeta(y)e^{-rh_\theta(y)}\,\mathrm dy
\asymp_{\alpha,n}
A_\zeta(y_\theta)y_\theta^{3/2}
e^{-rJ_\alpha(\theta)},
\label{eq:z-uniform-laplace-scaled}
\end{equation}
whenever $\lambda=r/(4y_\theta)\ge\rho/2$.
\end{lemma}

\begin{proof}
First observe that, for every $a\ge0$ and $v>0$,
\[
\min\{1,v\}
\le
\frac{1+av}{1+a}
\le
\max\{1,v\}.
\]
Applying this inequality to the two rational factors in
$R_{\zeta,\theta}$ yields, for some exponent $M=M(\alpha,n)>0$,
\[
R_{\zeta,\theta}(v)
\le
C_{\alpha,n}\bigl(v^M+v^{-M}\bigr),
\qquad v>0,
\]
uniformly in $\zeta\ge1$ and $\theta>0$. Moreover,
\[
R_{\zeta,\theta}(v)\asymp_{\alpha,n}1
\]
uniformly when $v$ belongs to a fixed neighborhood of $1$.

By \eqref{eq:z-normalized-phase}, the family $\Psi_\delta$ has a common minimum at $v=1$. The formula
\[
H_\delta''(v)
=
\frac{2}{v^3}
+(1-\delta)(q-1)v^{q-2}
\]
implies a uniform quadratic bound in a neighborhood of $1$. Outside that neighborhood, using in turn the term $1/v$, the term $\delta v$, and the term $(1-\delta)v^q/q$, there exist positive constants depending only on $\alpha$ such that
\[
\Psi_\delta(v)
\ge
c_\alpha
\begin{cases}
v^{-1},&0<v\le\frac12,\\
(v-1)^2,&\frac12\le v\le2,\\
v,&v\ge2.
\end{cases}
\]
In the central region one also has the upper bound
$\Psi_\delta(v)\le C_\alpha(v-1)^2$.

If $\lambda\ge1$, these estimates, together with the polynomial growth of
$R_{\zeta,\theta}$, show that the two tails are exponentially smaller than the contribution of a neighborhood of size $\lambda^{-1/2}$ around $v=1$. Integrating over this neighborhood gives the lower bound, while the upper bound follows from Gaussian localization in the central region and coercivity in the tails. Thus,
\[
\int_0^\infty
R_{\zeta,\theta}(v)
e^{-\lambda\Psi_{\delta_\theta}(v)}\,\mathrm dv
\asymp_{\alpha,n}\lambda^{-1/2},
\qquad \lambda\ge1.
\]
For $\rho/2\le\lambda\le1$, the right-hand side is comparable to one; the same amplitude and coercivity bounds provide uniform upper and lower estimates. This proves \eqref{eq:z-uniform-laplace}.

Finally, with $y=y_\theta v$ and
\eqref{eq:z-normalized-phase},
\[
\begin{aligned}
\int_0^\infty A_\zeta(y)e^{-rh_\theta(y)}\,\mathrm dy
&=
y_\theta A_\zeta(y_\theta)e^{-rJ_\alpha(\theta)}
\\
&\quad\times
\int_0^\infty
R_{\zeta,\theta}(v)
e^{-\lambda\Psi_{\delta_\theta}(v)}\,\mathrm dv,
\end{aligned}
\]
where $\lambda=r/(4y_\theta)$. Using \eqref{eq:z-uniform-laplace} and
$\lambda^{-1/2}=2y_\theta^{1/2}r^{-1/2}$ gives
\eqref{eq:z-uniform-laplace-scaled}.
\end{proof}

The critical-point equation \eqref{eq:z-saddle-equation} implies
\[
y_\theta\le\frac{1}{2\rho}.
\]
Therefore, in the exterior region $r\ge\max\{1,t^\alpha\}$,
\[
\lambda=\frac{r}{4y_\theta}
\ge
\frac{\rho r}{2}
\ge
\frac{\rho}{2},
\]
so the preceding lemma applies uniformly.

\begin{proposition}[Exterior two-sided estimate]\label{prop:z-exterior}
For every $t>0$ and
\[
r\ge\max\{1,t^\alpha\},
\]
one has
\begin{equation}
\begin{aligned}
Z_{\mathbb H}(t,r)
\asymp_{\alpha,n}{}&
t^{-\alpha}
\left(
\frac{1+y_{r/t}}{y_{r/t}}
\right)^{\frac{n-3}{2}}
\left(
1+\frac{r}{t^\alpha}y_{r/t}
\right)^\kappa
\\
&\times
\exp\!\left[
-r\left(
\rho+J_\alpha\!\left(\frac rt\right)
\right)
\right].
\end{aligned}
\label{eq:z-exterior-estimate}
\end{equation}
\end{proposition}

\begin{proof}
Apply \eqref{eq:z-uniform-laplace-scaled} to
\eqref{eq:z-exterior-reduction} with
\[
\zeta=\frac r{t^\alpha},
\qquad
\theta=\frac rt.
\]
Since
\[
A_\zeta(y_\theta)y_\theta^{3/2}
=
\left(
\frac{1+y_\theta}{y_\theta}
\right)^{\frac{n-3}{2}}
(1+\zeta y_\theta)^\kappa,
\]
we obtain \eqref{eq:z-exterior-estimate}.
\end{proof}

The rate $J_\alpha(r/t)$ is kept untruncated in
\eqref{eq:z-exterior-estimate}. This choice preserves the two-sided
equivalence uniformly throughout the exterior region; the expansions of
$J_\alpha$ in the limits $r/t\to0$ and $r/t\to\infty$ will be used only to
interpret subregimes, not to replace the global rate.

\subsection{Large time in the interior region}
\label{subsec:z-long-time-interior}

The large-time interior region is governed by the Green kernel
\eqref{eq:z-green-kernel}. We begin by recording its global size, which follows
directly from the Davies--Mandouvalos estimate.

\begin{proposition}[Green-kernel estimates]\label{prop:z-green-estimates}
For every $r>0$,
\begin{equation}
G_{\mathbb H}(r)
\asymp_n
\begin{cases}
1+|\log r|,
& n=2,\quad 0<r\le1,\\[1mm]
r^{2-n},
& n\ge3,\quad 0<r\le1,\\[1mm]
e^{-(n-1)r},
& r\ge1.
\end{cases}
\label{eq:z-green-estimates}
\end{equation}
\end{proposition}

\begin{proof}
Integrating \eqref{eq:heat-kernel-dm}, for $0<r\le1$ the tail $s\ge1$
is uniformly bounded from above, while the contribution of $s\in[1,2]$
is uniformly bounded from below. On $0<s\le1$, the geometric and spectral
factors other than $s^{-n/2}e^{-r^2/(4s)}$ are comparable to constants. Hence
\[
G_{\mathbb H}(r)
\asymp_n
1+
\int_0^1 s^{-n/2}e^{-r^2/(4s)}\,\mathrm ds.
\]
The change of variables $u=r^2/(4s)$ gives the logarithmic behavior when $n=2$ and
$r^{2-n}$ when $n\ge3$.

Now suppose $r\ge1$ and set $s=ru$. From
\eqref{eq:heat-kernel-dm} we obtain
\[
G_{\mathbb H}(r)
\asymp_n
r^{1/2}e^{-\rho r}
\int_0^\infty
u^{-n/2}(1+u)^{(n-3)/2}
\exp\!\left[-r\left(\rho^2u+\frac1{4u}\right)\right]
\,\mathrm du.
\]
The phase
\[
\varphi(u):=\rho^2u+\frac1{4u}
\]
is strictly convex and has a unique nondegenerate minimum at
$u_0=(2\rho)^{-1}$, with $\varphi(u_0)=\rho$. Laplace localization in
a fixed neighborhood of $u_0$, together with the coercivity of $\varphi$ at both
endpoints, gives
\[
\int_0^\infty
u^{-n/2}(1+u)^{(n-3)/2}e^{-r\varphi(u)}\,\mathrm du
\asymp_n r^{-1/2}e^{-\rho r}.
\]
Thus
\[
G_{\mathbb H}(r)\asymp_n e^{-2\rho r}=e^{-(n-1)r}.
\]
\end{proof}

\begin{lemma}[Uniform truncation of the Green kernel]
\label{lem:z-green-truncation}
There exists $c_n>0$ such that, for every $T\ge1$ and $0<r\le T$,
\begin{equation}
c_nG_{\mathbb H}(r)
\le
\int_0^{2T}p_s(r)\,\mathrm ds
\le
G_{\mathbb H}(r).
\label{eq:z-green-truncation}
\end{equation}
\end{lemma}

\begin{proof}
The upper bound is immediate from \eqref{eq:z-green-kernel}. For the lower
bound, we distinguish three cases.

If $0<r\le1$ and $n\ge3$, then
$[r^2,2r^2]\subset(0,2T]$. On this interval, by
\eqref{eq:heat-kernel-dm},
\[
p_s(r)\asymp_n s^{-n/2},
\]
and therefore
\[
\int_0^{2T}p_s(r)\,\mathrm ds
\ge
\int_{r^2}^{2r^2}p_s(r)\,\mathrm ds
\gtrsim_n r^{2-n}
\asymp_n G_{\mathbb H}(r).
\]

If $0<r\le1$ and $n=2$, for $r\le1/2$ we use
$s\in[r^2,1]$ and obtain
\[
\int_{r^2}^1p_s(r)\,\mathrm ds
\gtrsim 1+|\log r|.
\]
For $1/2\le r\le1$, the integral over $s\in[1,2]$ is uniformly
bounded from below. In both cases, \cref{prop:z-green-estimates} gives the
required bound.

Finally, let $1\le r\le T$ and define
\[
s_0:=\frac{r}{2\rho}=\frac{r}{n-1},
\qquad
I_r:=[s_0,s_0+\sqrt r].
\]
Since $n\ge2$ and $r\ge1$, one has $I_r\subset(0,2r]\subset(0,2T]$.
On $I_r$, $s\asymp_n r$, and the algebraic prefactor in
\eqref{eq:heat-kernel-dm} is comparable to $r^{-1/2}$. Moreover,
\[
\rho^2s+\frac{r^2}{4s}-\rho r
=
\frac{(2\rho s-r)^2}{4s},
\]
which remains bounded on $I_r$. Consequently,
\[
p_s(r)\gtrsim_n r^{-1/2}e^{-2\rho r},
\qquad s\in I_r.
\]
Integrating over an interval of length $\sqrt r$ and using
\cref{prop:z-green-estimates},
\[
\int_0^{2T}p_s(r)\,\mathrm ds
\gtrsim_n e^{-2\rho r}
\asymp_n G_{\mathbb H}(r).
\]
\end{proof}

\begin{proposition}[Large-time interior estimate]
\label{prop:z-long-interior}
For $t\ge1$ and $0<r\le t^\alpha$,
\begin{equation}
Z_{\mathbb H}(t,r)
\asymp_{\alpha,n}
t^{-\alpha}G_{\mathbb H}(r).
\label{eq:z-long-interior}
\end{equation}
\end{proposition}

\begin{proof}
By \eqref{eq:subordination-z},
\[
Z_{\mathbb H}(t,r)
=
t^{-\alpha}\int_0^\infty
p_s(r)M_\alpha(st^{-\alpha})\,\mathrm ds.
\]
The function $M_\alpha$ is bounded on $[0,\infty)$, so
\[
Z_{\mathbb H}(t,r)
\lesssim_\alpha
t^{-\alpha}G_{\mathbb H}(r).
\]
By positivity and continuity,
\[
c_\alpha^{(2)}:=\min_{0\le z\le2}M_\alpha(z)>0.
\]
Taking $T=t^\alpha\ge1$ and using $r\le T$,
\cref{lem:z-green-truncation} gives
\[
\begin{aligned}
Z_{\mathbb H}(t,r)
&\ge
c_\alpha^{(2)}t^{-\alpha}
\int_0^{2t^\alpha}p_s(r)\,\mathrm ds
\\
&\gtrsim_{\alpha,n}
t^{-\alpha}G_{\mathbb H}(r).
\end{aligned}
\]
\end{proof}

\subsection{Matching of the regimes and the global theorem}
\label{subsec:z-global}

To verify the matching between the local and exterior regimes, we will need
the following consequence of the critical-point equation.

\begin{lemma}[Critical point for $\theta\ge1$]\label{lem:z-saddle-large-theta}
For $\theta\ge1$,
\begin{equation}
y_\theta
\asymp_{\alpha,n}
\theta^{-\alpha/(2-\alpha)},
\label{eq:z-saddle-large-theta}
\end{equation}
and
\begin{equation}
0\le
J_\alpha(\theta)
-D_\alpha\theta^{\alpha/(2-\alpha)}
\le C_{\alpha,n}.
\label{eq:z-rate-large-theta}
\end{equation}
\end{lemma}

\begin{proof}
Write
\[
\beta:=\frac{\alpha}{2-\alpha}=\frac{m}{m+2}
\]
and set $y_\theta=\theta^{-\beta}x_\theta$ in
\eqref{eq:z-saddle-equation}. This gives
\[
\frac14
=
\rho^2\theta^{-2\beta}x_\theta^2
+\alpha^m x_\theta^{m+2}.
\]
Since $\theta\ge1$, this identity implies positive upper and lower
bounds for $x_\theta$ depending only on $(\alpha,n)$, which proves
\eqref{eq:z-saddle-large-theta}.

Let
\[
h_\theta^0(y):=\frac1{4y}+b_\alpha\theta^m y^q.
\]
Its minimum is
\[
\min_{y>0}h_\theta^0(y)
=
D_\alpha\theta^{\alpha/(2-\alpha)}.
\]
Since $h_\theta=h_\theta^0+\rho^2y$, the lower bound in
\eqref{eq:z-rate-large-theta} follows. If $y_\theta^0$ is the minimizer of
$h_\theta^0$, then
$y_\theta^0\asymp_\alpha\theta^{-\alpha/(2-\alpha)}\le C_\alpha$ for
$\theta\ge1$, and hence
\[
J_\alpha(\theta)
\le
h_\theta(y_\theta^0)
=
D_\alpha\theta^{\alpha/(2-\alpha)}+\rho^2y_\theta^0
\le
D_\alpha\theta^{\alpha/(2-\alpha)}+C_{\alpha,n}.
\]
\end{proof}

\begin{proposition}[Matching of the three regions]\label{prop:z-matching}
The expressions in \cref{prop:z-short-local,prop:z-exterior,prop:z-long-interior}
are uniformly comparable, with constants depending only on
$(\alpha,n)$, on the common interfaces $t=1$, $r=1$, and $r=t^\alpha$.
\end{proposition}

\begin{proof}
At $t=1$ and $0<r\le1$, \eqref{eq:z-local-near} and
\eqref{eq:z-green-estimates} show that
\[
F_{\alpha,n}(r^2)\asymp_{\alpha,n}G_{\mathbb H}(r).
\]
Thus the local and interior expressions are uniformly comparable. For
$r\ge1$, the exterior description is the same on both sides of $t=1$.

Now consider $r=t^\alpha$ with $t\ge1$. Then
\[
\theta=t^{\alpha-1}\in(0,1],
\qquad
\frac{r}{t^\alpha}=1.
\]
Equation \eqref{eq:z-saddle-equation} implies
$y_\theta\asymp_{\alpha,n}1$ uniformly for $0<\theta\le1$. Moreover,
since
\[
\min_{y>0}\left(\rho^2y+\frac1{4y}\right)=\rho,
\]
one has $J_\alpha(\theta)\ge\rho$, while evaluating
$h_\theta$ at $(2\rho)^{-1}$ gives
\[
0\le J_\alpha(\theta)-\rho\le C_{\alpha,n}\theta^m.
\]
On this interface,
\[
r\theta^m
=
t^\alpha(t^{\alpha-1})^{\alpha/(1-\alpha)}
=1.
\]
Therefore the exponential factor in \eqref{eq:z-exterior-estimate} is
comparable to $e^{-2\rho r}=e^{-(n-1)r}$, and its algebraic amplitude is
comparable to $t^{-\alpha}$. By \eqref{eq:z-green-estimates}, this agrees
with $t^{-\alpha}G_{\mathbb H}(r)$.

Finally, let $r=1$ and $0<t\le1$. Then
\[
\theta=t^{-1}\ge1,
\qquad
\Omega=t^{-\alpha}=\theta^\alpha.
\]
By \cref{lem:z-saddle-large-theta},
\[
y_\theta\asymp_{\alpha,n}\Omega^{-1/(2-\alpha)}
\]
and
\[
\left|
J_\alpha(\theta)
-
D_\alpha\Omega^{1/(2-\alpha)}
\right|
\le C_{\alpha,n}.
\]
Substituting these relations into \eqref{eq:z-exterior-estimate}, its
prefactor is comparable to
\[
t^{-n\alpha/[2(2-\alpha)]}.
\]
The local formula \eqref{eq:z-local-offdiag}, evaluated at
$\Omega=t^{-\alpha}$, has the same prefactor and the same exponential rate
up to a uniformly bounded multiplicative factor. This proves the
matching at $r=1$.
\end{proof}

Define the global profile
\begin{equation}
\mathcal Z_{\alpha,n}(t,r)
:=
\begin{cases}t^{-n\alpha/2}
F_{\alpha,n}\!\left(r^2/t^\alpha\right),
& 0<t\le1,\quad 0<r\le1,
\\[2mm]
t^{-\alpha}G_{\mathbb H}(r),
& t>1,\quad 0<r\le t^\alpha,
\\[2mm]
\begin{aligned}
&t^{-\alpha}
\left(
\dfrac{1+y_{r/t}}{y_{r/t}}
\right)^{\frac{n-3}{2}}
\left(
1+\dfrac r{t^\alpha}y_{r/t}
\right)^\kappa
\\[-1mm]
&\qquad\times
\exp\!\left[-r\left(\rho+J_\alpha(r/t)\right)\right]
\end{aligned},
& r>\max\{1,t^\alpha\}.
\end{cases}
\label{eq:z-global-profile}
\end{equation}
This choice makes the three cases disjoint. By \cref{prop:z-matching},
the alternative formulas available on the interfaces are equivalent with
uniform constants.

\begin{theorem}[Global two-sided estimate for $Z_{\mathbb H}$]
\label{thm:z-global}
Let $n\ge2$ and $0<\alpha<1$. There exist constants
$0<c_{\alpha,n}\le C_{\alpha,n}<\infty$ such that, for every $t>0$ and
$r>0$,
\begin{equation}
c_{\alpha,n}\mathcal Z_{\alpha,n}(t,r)
\le
Z_{\mathbb H}(t,r)
\le
C_{\alpha,n}\mathcal Z_{\alpha,n}(t,r),
\label{eq:z-global-estimate}
\end{equation}
where $\mathcal Z_{\alpha,n}$ is defined in
\eqref{eq:z-global-profile}. In particular, the constants are uniform over
all $(t,r)\in(0,\infty)^2$.
\end{theorem}

\begin{proof}
If $0<t\le1$ and $0<r\le1$, the claim is
\cref{prop:z-short-local}. If $t>1$ and $0<r\le t^\alpha$, it is
\cref{prop:z-long-interior}. In the region
$r>\max\{1,t^\alpha\}$, it is \cref{prop:z-exterior}. These three domains
cover $(0,\infty)^2$, and \cref{prop:z-matching} guarantees that the formulas
are compatible on their common interfaces. Therefore the constants can
be chosen to depend only on $(\alpha,n)$ throughout the whole domain.
\end{proof}

\section{Two-sided estimates for \texorpdfstring{$Y_{\mathbb H}$}{Y H}}
\label{sec:y-estimates}

In this section we study the inhomogeneous kernel $Y_{\mathbb H}(t,r)$. The subordination representation \eqref{eq:subordination-y} shows that, relative to $Z_{\mathbb H}$, an additional factor $s$ appears inside the integral. This modification does not change the phase in the master reduction, but it changes the algebraic powers, shifts the local dimensional threshold, and replaces the Green kernel by the potential kernel associated with $L^{-2}$. Every constant implicit in $\asymp_{\alpha,n}$ depends only on $(\alpha,n)$, unless explicitly stated otherwise.

\subsection{Master reduction}
\label{subsec:y-master-reduction}

\begin{proposition}[Master reduction for $Y_{\mathbb H}$]\label{prop:y-master-reduction}
For every $t>0$ and $r>0$,
\begin{equation}
\begin{aligned}
Y_{\mathbb H}(t,r)
\asymp_{\alpha,n}{}&
 t^{-\alpha-1}(1+r)e^{-\rho r}
\\
&\times
\int_0^\infty
s^{1-n/2}(1+r+s)^{(n-3)/2}
\left(1+\frac{s}{t^\alpha}\right)^\kappa
\exp\!\left(-\Phi_{t,r}(s)\right)\,\mathrm ds,
\end{aligned}
\label{eq:y-master-reduction}
\end{equation}
where $\Phi_{t,r}$ is the same phase defined in \eqref{eq:z-master-phase}.
\end{proposition}

\begin{proof}
By \eqref{eq:subordination-y},
\[
Y_{\mathbb H}(t,r)
=
\alpha t^{-\alpha-1}
\int_0^\infty
s\,p_s(r)M_\alpha(st^{-\alpha})\,\mathrm ds.
\]
Applying \eqref{eq:heat-kernel-dm} and \eqref{eq:m-wright-global} yields \eqref{eq:y-master-reduction}. The additional factor $s$ affects only the amplitude; the three exponential terms are still grouped into $\Phi_{t,r}$.
\end{proof}

\subsection{Governing profiles}
\label{subsec:y-governing-profiles}

In the short-time local region we retain
\[
\Omega:=\frac{r^2}{t^\alpha}
\]
and define, for $\Omega>0$,
\begin{equation}
F^Y_{\alpha,n}(\Omega)
:=
\int_0^\infty
\nu^{1-n/2}(1+\nu)^\kappa
\exp\!\left(
-\frac{\Omega}{4\nu}-b_\alpha\nu^q
\right)\,\mathrm d\nu.
\label{eq:y-local-profile}
\end{equation}
The only difference from \eqref{eq:z-local-profile} is the additional factor $\nu$ inherited from the subordination formula for $Y_{\mathbb H}$.

At large times in the interior region, the potential kernel associated with
$L^{-2}$,
\begin{equation}
G_{\mathbb H}^{(2)}(r)
:=
\int_0^\infty s\,p_s(r)\,\mathrm ds,
\qquad r>0,
\label{eq:y-second-potential}
\end{equation}
appears. Indeed, $G_{\mathbb H}^{(2)}$ is the kernel of $L^{-2}$, since
\[
\int_0^\infty s e^{-sL}\,\mathrm ds=L^{-2},
\qquad L=-\Delta_{\mathbb H^n}.
\]
These two objects will play for $Y_{\mathbb H}$ the same roles that $F_{\alpha,n}$ and $G_{\mathbb H}$ play for $Z_{\mathbb H}$.

\subsection{Short time in the local region}
\label{subsec:y-short-time-local}

Throughout this subsection we assume
\[
0<t\le1,
\qquad
0<r\le1,
\]
and write
\[
\Omega=\frac{r^2}{t^\alpha},
\qquad
K=t^{-\alpha}\ge1.
\]
Then $0<\Omega\le K$. The change of variables $s=t^\alpha\nu$ in
\eqref{eq:y-master-reduction} gives
\begin{equation}
Y_{\mathbb H}(t,r)
\asymp_{\alpha,n}
t^{\alpha-1-n\alpha/2}
\int_0^\infty
W_{t,r}(\nu)\,f^Y_\Omega(\nu)\,\mathrm d\nu,
\label{eq:y-short-rescaled}
\end{equation}
where $W_{t,r}$ is the same geometric factor as in
\eqref{eq:z-short-rescaled}, and
\[
f^Y_\Omega(\nu)
:=
\nu^{1-n/2}(1+\nu)^\kappa
\exp\!\left(
-\frac{\Omega}{4\nu}-b_\alpha\nu^q
\right).
\]
In particular,
\[
0<W_{t,r}(\nu)\le C_n,
\qquad \nu>0,
\]
and
\[
c_n\le W_{t,r}(\nu)\le C_n,
\qquad 0<\nu\le K.
\]
The comparison between the truncated and full profiles required for $Y_{\mathbb H}$ is the same as for $Z_{\mathbb H}$, except for one additional power of $\nu$.

\begin{lemma}[Uniform truncation of the local profile of $Y_{\mathbb H}$]
\label{lem:y-local-truncation}
Let
\[
I^Y_K(\Omega)
:=
\int_0^K f^Y_\Omega(\nu)\,\mathrm d\nu.
\]
There exist constants $0<c_{\alpha,n}\le C_{\alpha,n}<\infty$ such that
\begin{equation}
c_{\alpha,n}F^Y_{\alpha,n}(\Omega)
\le
I^Y_K(\Omega)
\le
F^Y_{\alpha,n}(\Omega)
\le
C_{\alpha,n}I^Y_K(\Omega)
\label{eq:y-local-truncation}
\end{equation}
for every $K\ge1$ and $0<\Omega\le K$.
\end{lemma}

\begin{proof}
The argument is the same as in \cref{lem:z-local-truncation}, the only difference being that $f^Y_\Omega(\nu)=\nu f_\Omega(\nu)$. The additional factor $\nu$ changes only algebraic powers. Indeed, for
\[
T^Y_K(\Omega):=\int_K^\infty f^Y_\Omega(\nu)\,\mathrm d\nu
\]
one obtains, with exponents depending only on $(\alpha,n)$,
\[
T^Y_K(\Omega)
\le C_{\alpha,n}K^{B_1}e^{-c_\alpha K^q}.
\]
If $1\le K\le K_0$, or if $K\ge K_0$ and $0<\Omega\le1$, the interval
$\nu\in[1/2,1]$ gives the same uniform lower bound as for $Z_{\mathbb H}$.
If $K\ge K_0$ and $1\le\Omega\le K$, the same critical point
$\nu_\Omega\asymp_\alpha\Omega^{(1-\alpha)/(2-\alpha)}$ satisfies
$2\nu_\Omega\le K$ for $K_0$ large enough, and integration over
$[\nu_\Omega,2\nu_\Omega]$ yields
\[
I^Y_K(\Omega)
\ge c\,K^{-A}\exp\!\left(-C K^{1/(2-\alpha)}\right).
\]
Since $q>1/(2-\alpha)$, this lower bound uniformly absorbs the preceding tail. This proves \eqref{eq:y-local-truncation}.
\end{proof}

\begin{proposition}[Short-time local regime for $Y_{\mathbb H}$]
\label{prop:y-short-local}
For $0<t\le1$ and $0<r\le1$,
\begin{equation}
Y_{\mathbb H}(t,r)
\asymp_{\alpha,n}
t^{\alpha-1-n\alpha/2}
F^Y_{\alpha,n}\!\left(\frac{r^2}{t^\alpha}\right).
\label{eq:y-short-local}
\end{equation}
\end{proposition}

\begin{proof}
The upper bound follows from \eqref{eq:y-short-rescaled} using
$W_{t,r}\le C_n$. For the lower bound, we restrict to $0<\nu\le K$,
use $W_{t,r}\ge c_n$, and apply \cref{lem:y-local-truncation}.
\end{proof}

The additional factor $\nu$ shifts the dimensional threshold of the local profile by two dimensions.

\begin{proposition}[Two-sided asymptotics of the local profile of $Y_{\mathbb H}$]
\label{prop:y-local-profile-asymptotics}
For $0<\Omega\le1$,
\begin{equation}
F^Y_{\alpha,n}(\Omega)
\asymp_{\alpha,n}
\begin{cases}
1,
&2\le n\le3,\\[1mm]
1+|\log\Omega|,
&n=4,\\[1mm]
\Omega^{2-n/2},
&n\ge5.
\end{cases}
\label{eq:y-local-near}
\end{equation}
For $\Omega\ge1$,
\begin{equation}
F^Y_{\alpha,n}(\Omega)
\asymp_{\alpha,n}
\Omega^{-\frac{(n-2)(1-\alpha)}{2(2-\alpha)}}
\exp\!\left(
-D_\alpha\Omega^{1/(2-\alpha)}
\right),
\label{eq:y-local-offdiag}
\end{equation}
where $D_\alpha$ is the constant in \eqref{eq:z-local-rate-constant}.
\end{proposition}

\begin{proof}
If $0<\Omega\le1$, the contribution of $\nu\ge1$ to
\eqref{eq:y-local-profile} is comparable to a constant. On $0<\nu\le1$,
the factors $(1+\nu)^\kappa$ and $e^{-b_\alpha\nu^q}$ are comparable to
one, and therefore
\[
F^Y_{\alpha,n}(\Omega)
\asymp_{\alpha,n}
1+
\int_0^1
\nu^{1-n/2}e^{-\Omega/(4\nu)}\,\mathrm d\nu.
\]
With $u=\Omega/(4\nu)$,
\[
\int_0^1
\nu^{1-n/2}e^{-\Omega/(4\nu)}\,\mathrm d\nu
=
4^{n/2-2}\Omega^{2-n/2}
\int_{\Omega/4}^{\infty}u^{n/2-3}e^{-u}\,\mathrm du.
\]
The incomplete integral is of order $\Omega^{n/2-2}$ for $n<4$,
of order $1+|\log\Omega|$ for $n=4$, and of order one for $n>4$.
This gives \eqref{eq:y-local-near}.

For $\Omega\ge1$, we use the same critical point $\nu_\Omega$ and the same
parameter $\lambda_\Omega$ as in the proof of
\cref{prop:z-local-profile-asymptotics}. After $\nu=\nu_\Omega v$,
the normalized phase is identical, while the amplitude gains a factor
$\nu_\Omega v$. The same uniform Laplace argument yields
\[
F^Y_{\alpha,n}(\Omega)
\asymp_{\alpha,n}
\nu_\Omega^{2-n/2}(1+\nu_\Omega)^\kappa
\lambda_\Omega^{-1/2}
\exp\!\left(-D_\alpha\Omega^{1/(2-\alpha)}\right).
\]
Since $\nu_\Omega\asymp_\alpha
\Omega^{(1-\alpha)/(2-\alpha)}$, the algebraic power reduces to that in
\eqref{eq:y-local-offdiag}.
\end{proof}

Combining \cref{prop:y-short-local,prop:y-local-profile-asymptotics},
for $0<t\le1$ and $0<r\le t^{\alpha/2}$ we obtain
\[
Y_{\mathbb H}(t,r)
\asymp_{\alpha,n}
\begin{cases}t^{\alpha-1-n\alpha/2},
&2\le n\le3,\\[1mm]
t^{-\alpha-1}
\left(
1+\left|\log\dfrac{r^2}{t^\alpha}\right|
\right),
&n=4,\\[2mm]
t^{-\alpha-1}r^{4-n},
&n\ge5,
\end{cases}
\]
whereas, uniformly for
$0<t\le1$ and $t^{\alpha/2}\le r\le1$,
\[
\begin{aligned}
Y_{\mathbb H}(t,r)
\asymp_{\alpha,n}{}&
t^{\alpha-1-n\alpha/2}
\left(\frac{r^2}{t^\alpha}\right)^{
-\frac{(n-2)(1-\alpha)}{2(2-\alpha)}
}
\\
&\times
\exp\!\left[
-D_\alpha
\left(\frac{r^2}{t^\alpha}\right)^{1/(2-\alpha)}
\right].
\end{aligned}
\]
In particular, the local dimensional threshold for $Y_{\mathbb H}$ is $n=4$. The same threshold appears in the estimates for the Euclidean inhomogeneous kernel of Eidelman--Kochubei, where $Y$ is the Riemann--Liouville derivative of the kernel associated with the initial datum \cite{EidelmanKochubei2004}.

\subsection{Large time in the interior region}
\label{subsec:y-long-time-interior}

The large-time interior region is governed by the potential
$G_{\mathbb H}^{(2)}$. We first establish its two-sided estimates and a
uniform truncated version, which will be the direct input in the
subordination formula for $Y_{\mathbb H}$.

\begin{proposition}[Estimates for the $L^{-2}$ potential kernel]
\label{prop:y-second-potential-estimates}
For every $r>0$,
\begin{equation}
G_{\mathbb H}^{(2)}(r)
\asymp_n
\begin{cases}
1,
&2\le n\le3,\quad 0<r\le1,\\[1mm]
1+|\log r|,
&n=4,\quad 0<r\le1,\\[1mm]
r^{4-n},
&n\ge5,\quad 0<r\le1,\\[1mm]
r e^{-(n-1)r},
&r\ge1.
\end{cases}
\label{eq:y-second-potential-estimates}
\end{equation}
\end{proposition}

\begin{proof}
First suppose $0<r\le1$. By \eqref{eq:heat-kernel-dm}, for
$0<s\le1$ the geometric and spectral factors other than
$s^{-n/2}e^{-r^2/(4s)}$ are uniformly comparable to one. Moreover, the
contribution of $s\ge1$ to \eqref{eq:y-second-potential} is comparable to
a positive constant, uniformly for $0<r\le1$. Hence
\[
G_{\mathbb H}^{(2)}(r)
\asymp_n
1+
\int_0^1
s^{1-n/2}e^{-r^2/(4s)}\,\mathrm ds.
\]
With the change of variables $u=r^2/(4s)$,
\[
\int_0^1
s^{1-n/2}e^{-r^2/(4s)}\,\mathrm ds
=
4^{n/2-2}r^{4-n}
\int_{r^2/4}^{\infty}u^{n/2-3}e^{-u}\,\mathrm du.
\]
For $2\le n\le3$, the incomplete integral is of order $r^{n-4}$ and the
product with the prefactor $r^{4-n}$ remains bounded above and below. For
$n=4$, the growth $1+|\log r|$ appears, whereas for $n\ge5$ the
incomplete integral is comparable to a positive constant and the factor
$r^{4-n}$ remains. This proves the first three lines of
\eqref{eq:y-second-potential-estimates}.

Now let $r\ge1$ and set $s=ry$. The Davies--Mandouvalos estimate gives
\[
G_{\mathbb H}^{(2)}(r)
\asymp_n
r^{3/2}e^{-\rho r}
\int_0^\infty
 y^{1-n/2}(1+y)^{(n-3)/2}
\exp\!\left[-r\left(\rho^2y+\frac1{4y}\right)\right]\,\mathrm dy.
\]
The phase $\rho^2y+1/(4y)$ has a unique nondegenerate minimum at
$y=(2\rho)^{-1}$, with value $\rho$. The uniform Laplace argument already
used in \cref{lem:z-uniform-laplace}, here in the limiting case
$\delta=1$, yields a factor $r^{-1/2}$. Consequently,
\[
G_{\mathbb H}^{(2)}(r)
\asymp_n
r e^{-2\rho r}
=
r e^{-(n-1)r}.
\]
\end{proof}

\begin{lemma}[Uniform weighted truncation]
\label{lem:y-second-potential-truncation}
There exists $c_n>0$ such that, for every $T\ge1$ and $0<r\le T$,
\begin{equation}
c_n G_{\mathbb H}^{(2)}(r)
\le
\int_0^{2T}s\,p_s(r)\,\mathrm ds
\le
G_{\mathbb H}^{(2)}(r).
\label{eq:y-second-potential-truncation}
\end{equation}
\end{lemma}

\begin{proof}
The upper bound is immediate from \eqref{eq:y-second-potential}. For the
lower bound we distinguish the same regimes as in
\eqref{eq:y-second-potential-estimates}.

If $0<r\le1$ and $n\ge5$, we integrate over
$s\in[r^2,2r^2]\subset(0,2T]$. On this interval
$p_s(r)\asymp_n s^{-n/2}$, and hence
\[
\int_{r^2}^{2r^2}s\,p_s(r)\,\mathrm ds
\asymp_n r^{4-n}.
\]
If $n=4$, integration over $[r^2,1]$, together with a contribution from
a fixed interval contained in $[1,2]$, gives
$1+|\log r|$. If $2\le n\le3$, the contribution of a fixed interval,
for instance $s\in[1,2]$, is sufficient.

Finally, if $1\le r\le T$, set
\[
s_0:=\frac{r}{n-1},
\qquad
I_r:=[s_0,s_0+\sqrt r]\subset(0,2T].
\]
On $I_r$ one has $s\asymp_n r$ and, as in the proof of
\cref{lem:z-green-truncation},
\[
p_s(r)\gtrsim_n r^{-1/2}e^{-(n-1)r}.
\]
Since $|I_r|=\sqrt r$,
\[
\int_{I_r}s\,p_s(r)\,\mathrm ds
\gtrsim_n
r e^{-(n-1)r}
\asymp_n G_{\mathbb H}^{(2)}(r).
\]
This completes the proof.
\end{proof}

\begin{proposition}[Large-time interior regime for $Y_{\mathbb H}$]
\label{prop:y-long-interior}
For $t\ge1$ and $0<r\le t^\alpha$,
\begin{equation}
Y_{\mathbb H}(t,r)
\asymp_{\alpha,n}
t^{-\alpha-1}G_{\mathbb H}^{(2)}(r).
\label{eq:y-long-interior}
\end{equation}
\end{proposition}

\begin{proof}
By \eqref{eq:subordination-y},
\[
Y_{\mathbb H}(t,r)
=
\alpha t^{-\alpha-1}
\int_0^\infty
s\,p_s(r)M_\alpha(st^{-\alpha})\,\mathrm ds.
\]
Since $M_\alpha$ is bounded on $[0,\infty)$,
\[
Y_{\mathbb H}(t,r)
\lesssim_{\alpha,n}
t^{-\alpha-1}G_{\mathbb H}^{(2)}(r).
\]
By positivity and continuity of $M_\alpha$,
\[
c_{\alpha,2}^{M}:=
\min_{0\le z\le2}M_\alpha(z)>0.
\]
Taking $T=t^\alpha\ge1$ and using
\cref{lem:y-second-potential-truncation}, for $r\le t^\alpha$ we obtain
\[
\begin{aligned}
Y_{\mathbb H}(t,r)
&\ge
\alpha c_{\alpha,2}^{M}t^{-\alpha-1}
\int_0^{2t^\alpha}s\,p_s(r)\,\mathrm ds
\\
&\gtrsim_{\alpha,n}
t^{-\alpha-1}G_{\mathbb H}^{(2)}(r).
\end{aligned}
\]
This proves \eqref{eq:y-long-interior}.
\end{proof}

\subsection{Transfer of the exterior critical point}
\label{subsec:y-exterior-transfer}

Consider
\[
r\ge\max\{1,t^\alpha\},
\qquad
\zeta:=\frac{r}{t^\alpha}\ge1,
\qquad
\theta:=\frac rt.
\]
Making the change of variables $s=ry$ in \eqref{eq:y-master-reduction} and using the same algebraic comparisons as in \eqref{eq:z-exterior-reduction}, we obtain
\begin{equation}
Y_{\mathbb H}(t,r)
\asymp_{\alpha,n}
 t^{-\alpha-1}r^{3/2}e^{-\rho r}
\int_0^\infty
B_\zeta(y)e^{-r h_\theta(y)}\,\mathrm dy,
\label{eq:y-exterior-reduction}
\end{equation}
where $h_\theta$ is the phase in \eqref{eq:z-exterior-phase} and
\begin{equation}
B_\zeta(y)
:=
y^{1-n/2}(1+y)^{(n-3)/2}(1+\zeta y)^\kappa
=
yA_\zeta(y).
\label{eq:y-exterior-amplitude}
\end{equation}
Thus $Y_{\mathbb H}$ has the same minimizer $y_\theta$ and the same rate $J_\alpha(\theta)$ as $Z_{\mathbb H}$.

\begin{proposition}[Transfer of the exterior critical point]\label{prop:y-exterior-transfer}
For every $t>0$ and $r\ge\max\{1,t^\alpha\}$,
\begin{equation}
\begin{aligned}
Y_{\mathbb H}(t,r)
\asymp_{\alpha,n}{}&
 t^{-\alpha-1}r\,y_{r/t}
\left(
\frac{1+y_{r/t}}{y_{r/t}}
\right)^{\frac{n-3}{2}}
\left(
1+\frac r{t^\alpha}y_{r/t}
\right)^\kappa
\\
&\times
\exp\!\left[
-r\left(
\rho+J_\alpha\!\left(\frac rt\right)
\right)
\right].
\end{aligned}
\label{eq:y-exterior-estimate}
\end{equation}
In particular, uniformly throughout the exterior region,
\begin{equation}
Y_{\mathbb H}(t,r)
\asymp_{\alpha,n}
\frac{r\,y_{r/t}}{t}\,
Z_{\mathbb H}(t,r).
\label{eq:y-z-exterior-relation}
\end{equation}
\end{proposition}

\begin{proof}
By \eqref{eq:y-exterior-amplitude},
\[
\frac{B_\zeta(y_\theta v)}{B_\zeta(y_\theta)}
=
v\,
\frac{A_\zeta(y_\theta v)}{A_\zeta(y_\theta)}.
\]
The additional factor $v$ preserves the uniform polynomial control used in \cref{lem:z-uniform-laplace}. The normalized phase, its coercivity, and the parameter $\lambda=r/(4y_\theta)$ are the same. Therefore,
\[
r^{1/2}
\int_0^\infty
B_\zeta(y)e^{-rh_\theta(y)}\,\mathrm dy
\asymp_{\alpha,n}
B_\zeta(y_\theta)y_\theta^{3/2}e^{-rJ_\alpha(\theta)}.
\]
Since $B_\zeta(y_\theta)=y_\theta A_\zeta(y_\theta)$, substituting this estimate into \eqref{eq:y-exterior-reduction} gives \eqref{eq:y-exterior-estimate}. Finally, comparison with \eqref{eq:z-exterior-estimate} yields \eqref{eq:y-z-exterior-relation}.
\end{proof}

The two-sided relation \eqref{eq:y-z-exterior-relation} summarizes the difference between the two kernels in the exterior region: the exponential rate is the same, whereas the amplitude of $Y_{\mathbb H}$ acquires the factor $r y_{r/t}/t$ evaluated at the same critical point.

\subsection{Global theorem and matching}
\label{subsec:y-global}

The three estimates obtained above cover the whole domain $(t,r)\in(0,\infty)^2$.
Before gathering them into a single statement, we verify that the expressions
meeting at the common interfaces are uniformly comparable.

\begin{proposition}[Matching of the three regions]\label{prop:y-matching}
The expressions in
\cref{prop:y-short-local,prop:y-long-interior,prop:y-exterior-transfer}
are uniformly comparable, with constants depending only on
$(\alpha,n)$, on the common interfaces $t=1$, $r=1$, and $r=t^\alpha$.
\end{proposition}

\begin{proof}
At $t=1$ and $0<r\le1$, \eqref{eq:y-local-near} and
\eqref{eq:y-second-potential-estimates} give
\[
F^Y_{\alpha,n}(r^2)
\asymp_{\alpha,n}
G_{\mathbb H}^{(2)}(r).
\]
Indeed, both sides are of order one for $2\le n\le3$, of order
$1+|\log r|$ for $n=4$, and of order $r^{4-n}$ for $n\ge5$. Hence the
local and interior expressions are uniformly comparable. For $r\ge1$, the
same exterior formula applies on both sides of $t=1$.

Now consider $r=t^\alpha$ with $t\ge1$. Then
\[
\theta=t^{\alpha-1}\in(0,1],
\qquad
\frac{r}{t^\alpha}=1.
\]
As in the proof of \cref{prop:z-matching}, the critical-point equation gives
\[
y_\theta\asymp_{\alpha,n}1,
\qquad
0\le J_\alpha(\theta)-\rho
\le C_{\alpha,n}\theta^m,
\]
and moreover $r\theta^m=1$. Consequently, the exponential factor in
\eqref{eq:y-exterior-estimate} is comparable to
$e^{-2\rho r}=e^{-(n-1)r}$, whereas its algebraic amplitude is
comparable to
\[
t^{-\alpha-1}r=t^{-1}.
\]
By \eqref{eq:y-second-potential-estimates}, the interior expression at this
interface satisfies
\[
t^{-\alpha-1}G_{\mathbb H}^{(2)}(r)
\asymp_n
t^{-\alpha-1}r e^{-(n-1)r}
=
t^{-1}e^{-(n-1)r},
\]
which proves the second matching.

Finally, let $r=1$ and $0<t\le1$. Write
\[
\theta=t^{-1}\ge1,
\qquad
\Omega=t^{-\alpha}=\theta^\alpha.
\]
By \eqref{eq:y-local-offdiag} and \eqref{eq:z-local-offdiag}, the ratio between
the local profiles of $Y_{\mathbb H}$ and $Z_{\mathbb H}$ is comparable to
\[
t^{\alpha-1}
\Omega^{\frac{1-\alpha}{2-\alpha}}
=
t^{-\frac{2(1-\alpha)}{2-\alpha}}.
\]
On the other hand, \eqref{eq:y-z-exterior-relation} and
\eqref{eq:z-saddle-large-theta} give, for the ratio of the exterior profiles,
\[
\frac{y_{1/t}}{t}
\asymp_{\alpha,n}
t^{-1+\alpha/(2-\alpha)}
=
t^{-\frac{2(1-\alpha)}{2-\alpha}}.
\]
Since the local and exterior profiles of $Z_{\mathbb H}$ match uniformly
at $r=1$ by \cref{prop:z-matching}, the same holds for
$Y_{\mathbb H}$.
\end{proof}

Define the global profile
\begin{equation}
\mathcal Y_{\alpha,n}(t,r)
:=
\begin{cases}t^{\alpha-1-n\alpha/2}
F^Y_{\alpha,n}\!\left(r^2/t^\alpha\right),
& 0<t\le1,\quad 0<r\le1,
\\[2mm]
t^{-\alpha-1}G_{\mathbb H}^{(2)}(r),
& t>1,\quad 0<r\le t^\alpha,
\\[2mm]
\begin{aligned}
&t^{-\alpha-1}r\,y_{r/t}
\left(
\dfrac{1+y_{r/t}}{y_{r/t}}
\right)^{\frac{n-3}{2}}
\left(
1+\dfrac r{t^\alpha}y_{r/t}
\right)^\kappa
\\[-1mm]
&\qquad\times
\exp\!\left[-r\left(\rho+J_\alpha(r/t)\right)\right]
\end{aligned},
& r>\max\{1,t^\alpha\}.
\end{cases}
\label{eq:y-global-profile}
\end{equation}
The three cases are disjoint and cover the whole domain. By
\cref{prop:y-matching}, the alternative formulas available on the
interfaces are uniformly comparable.

\begin{theorem}[Global two-sided estimate for $Y_{\mathbb H}$]
\label{thm:y-global}
Let $n\ge2$ and $0<\alpha<1$. There exist constants
$0<c_{\alpha,n}\le C_{\alpha,n}<\infty$ such that, for every $t>0$ and
$r>0$,
\begin{equation}
c_{\alpha,n}\mathcal Y_{\alpha,n}(t,r)
\le
Y_{\mathbb H}(t,r)
\le
C_{\alpha,n}\mathcal Y_{\alpha,n}(t,r),
\label{eq:y-global-estimate}
\end{equation}
where $\mathcal Y_{\alpha,n}$ is defined in
\eqref{eq:y-global-profile}. In particular, the constants are uniform over
all $(t,r)\in(0,\infty)^2$.
\end{theorem}

\begin{proof}
If $0<t\le1$ and $0<r\le1$, the claim is
\cref{prop:y-short-local}. If $t>1$ and $0<r\le t^\alpha$, it is
\cref{prop:y-long-interior}. In the region $r>\max\{1,t^\alpha\}$, it is
\cref{prop:y-exterior-transfer}. These three domains cover
$(0,\infty)^2$, and \cref{prop:y-matching} guarantees uniform compatibility
on the common interfaces. Therefore the constants can be chosen to
depend only on $(\alpha,n)$ throughout the whole domain.
\end{proof}

\section{Geometric and probabilistic consequences}
\label{sec:properties}

The subordination representations, together with the geometry of the classical
hyperbolic heat flow, allow us to isolate properties that complement the
pointwise estimates above but do not logically depend on them. In particular,
subordination preserves the mass structure of $Z_{\mathbb H}$, yields a
natural normalization of $Y_{\mathbb H}$, and transfers the radial propagation
of hyperbolic heat to fractional time. This section collects these consequences
and makes explicit their relation to the geometry of $\mathbb H^n$.

\subsection{Mass and probabilistic normalization}
\label{subsec:properties-mass}

The mass identity \eqref{eq:kernel-masses} shows a structural difference
between the two kernels:
\[
\int_{\mathbb H^n} Z_{\mathbb H}(t,d(o,x))\,\mathrm d\mu(x)=1,
\qquad
\int_{\mathbb H^n} Y_{\mathbb H}(t,d(o,x))\,\mathrm d\mu(x)
=\frac{t^{\alpha-1}}{\Gamma(\alpha)}.
\]
Thus, $Z_{\mathbb H}(t)$ preserves mass, whereas $Y_{\mathbb H}(t)$ has
the mass prescribed by the temporal Duhamel kernel.

Define the normalized kernel
\begin{equation}
\widetilde Y_{\mathbb H}(t,r)
:=
\Gamma(\alpha)t^{1-\alpha}Y_{\mathbb H}(t,r).
\label{eq:y-normalized-kernel}
\end{equation}
Then $\widetilde Y_{\mathbb H}(t,\cdot)$ has mass one. Moreover, after the
change of variables $s=t^\alpha z$ in \eqref{eq:subordination-y},
\begin{equation}
\widetilde Y_{\mathbb H}(t)
=
\int_0^\infty
 e^{t^\alpha z\Delta_{\mathbb H^n}}
\,\Gamma(1+\alpha)zM_\alpha(z)\,\mathrm dz.
\label{eq:y-size-biased-subordination}
\end{equation}
By the moment formula for $M_\alpha$,
\[
\int_0^\infty \Gamma(1+\alpha)zM_\alpha(z)\,\mathrm dz=1.
\]
Therefore, the measure
\[
\nu_\alpha(\mathrm dz)
:=
\Gamma(1+\alpha)zM_\alpha(z)\,\mathrm dz
\]
is a probability measure. If $\mathcal{E}_1$ denotes a random variable with
density $M_\alpha(z)\,\mathrm dz$, then $\nu_\alpha$ is precisely its
size-biased law. In this sense, $Z_{\mathbb H}$ averages the heat semigroup
with respect to the inverse stable clock, whereas the normalization of
$Y_{\mathbb H}$ averages it with respect to the size-biased version of the
same clock.

\subsection{Inverse clock and radial propagation}
\label{subsec:properties-radial}

Let $(B_s)_{s\ge0}$ be hyperbolic Brownian motion with generator
$\Delta_{\mathbb H^n}$, started at a pole $o$, and let
$(\mathcal{E}_t)_{t\ge0}$ be the inverse $\alpha$--stable clock, independent
of $B$, with density $\ell_\alpha(t,s)$. We denote by $(X_t)_{t\ge0}$ the
time-fractional diffusion whose law at time $t$, started at $o$, has density
$x\mapsto Z_{\mathbb H}(t,d(o,x))$. Formula
\eqref{eq:subordination-z} is equivalent to the distributional representation,
standard in the probabilistic interpretation of fractional Cauchy problems by
inverse subordination
\cite{BaeumerMeerschaert2001,MeerschaertNaneVellaisamy2009}:
\begin{equation}
X_t\overset{d}=B_{\mathcal{E}_t}.
\label{eq:time-change-representation}
\end{equation}
Here $\overset d=$ denotes equality in distribution; below we use
$\Longrightarrow$ for convergence in distribution. We write
\[
R_s^{\rm cl}:=d(o,B_s),
\qquad
R_t:=d(o,X_t),
\]
and, for $p>0$,
\[
m_p(s):=\mathbb E[(R_s^{\rm cl})^p].
\]
By Tonelli's theorem and the self-similarity
$\mathcal{E}_t\overset d=t^\alpha \mathcal{E}_1$,
\begin{equation}
\mathbb E[R_t^p]
=
\int_0^\infty m_p(s)\ell_\alpha(t,s)\,\mathrm ds
=
\int_0^\infty m_p(t^\alpha z)M_\alpha(z)\,\mathrm dz.
\label{eq:radial-moment-subordination}
\end{equation}

For the classical heat flow, hyperbolic geometry produces ballistic radial
propagation. In particular, Vázquez shows that, at large times, the mass of
the heat kernel concentrates around $r=(n-1)s$ with width of order
$s^{1/2}$ \cite{Vazquez2022}. In the notation above,
\[
\frac{R_s^{\rm cl}}{s}\longrightarrow n-1
\qquad\text{in probability as }s\to\infty.
\]
Moreover, the Davies--Mandouvalos estimate implies, for every $p>0$,
\begin{equation*}
m_p(s)\le C_{p,n}\bigl(s^{p/2}+s^p\bigr),
\qquad s>0.
\end{equation*}
We include the argument because this domination will be used twice. For
$0<s\le1$, we split the radial integral into $r\le1$ and $r\ge1$. In the
first region, $\sinh r\asymp r$ and the Gaussian bound from
\eqref{eq:heat-kernel-dm}, followed by the change of variables
$r=\sqrt{s}\,u$, gives
$m_p(s)\lesssim_{p,n}s^{p/2}$; the tail $r\ge1$ is absorbed by the same
Gaussian factor. For $s\ge1$, in the relevant region $r\ge1$ we use
$\sinh r\asymp e^r$ and complete the square in the phase, obtaining a
Gaussian centered at $r=(n-1)s$ with width $\sqrt{s}$. The change of
variables $r=(n-1)s+\sqrt{s}\,u$ gives
$m_p(s)\lesssim_{p,n}s^p$, while the contribution of $r\le1$ is smaller.
This yields the global bound above. Applying it with an exponent slightly
larger than $p$ also gives the uniform integrability needed to transfer the
moments from the ballistic regime. In particular, for $t\ge1$,
\[
\frac{m_p(t^\alpha z)}{t^{\alpha p}}
\le C_{p,n}\left(z^{p/2}+z^p\right),
\]
and, for $0<t\le1$,
\[
\frac{m_2(t^\alpha z)}{t^\alpha}
\le C_n(z+z^2).
\]
Both right-hand sides are integrable with respect to
$M_\alpha(z)\,\mathrm dz$ by the moment formula for $M_\alpha$.

\begin{proposition}[Fractional radial law and moments]
\label{prop:radial-law-moments}
For $0<\alpha<1$,
\begin{equation}
\frac{R_t}{t^\alpha}
\Longrightarrow
(n-1)\mathcal{E}_1,
\qquad t\to\infty,
\label{eq:radial-law}
\end{equation}
and, for every $p>0$,
\begin{equation}
\mathbb E[R_t^p]
\sim
(n-1)^p
\frac{\Gamma(p+1)}{\Gamma(1+\alpha p)}
 t^{\alpha p},
\qquad t\to\infty.
\label{eq:radial-moments}
\end{equation}
\end{proposition}

\begin{proof}
To prove the one-time marginal statement, introduce an auxiliary coupling
$\widehat{\mathcal{E}}_t:=t^\alpha \mathcal{E}_1$, with
$\mathcal{E}_1$ independent of $B$. For every $t>0$,
$\widehat{\mathcal{E}}_t\overset d=\mathcal{E}_t$; we make no claim here
about equality of the full processes. Since $\mathcal{E}_1>0$ almost surely,
$\widehat{\mathcal{E}}_t\to\infty$ almost surely. The ballistic concentration
of the classical heat flow then implies
$R_{\widehat{\mathcal{E}}_t}^{\rm cl}/\widehat{\mathcal{E}}_t\to n-1$ in
probability. Multiplying by
$\widehat{\mathcal{E}}_t/t^\alpha=\mathcal{E}_1$ and using equality of the
one-time marginal laws yields \eqref{eq:radial-law}.

For the moments, divide \eqref{eq:radial-moment-subordination} by
$t^{\alpha p}$. For every $z>0$,
\[
\frac{m_p(t^\alpha z)}{t^{\alpha p}}
\longrightarrow
(n-1)^p z^p.
\]
The explicit domination above, valid for $t\ge1$, provides an integrable
function with respect to $M_\alpha(z)\,\mathrm dz$, and dominated convergence
gives
\[
\frac{\mathbb E[R_t^p]}{t^{\alpha p}}
\longrightarrow
(n-1)^p\int_0^\infty z^pM_\alpha(z)\,\mathrm dz.
\]
The moment formula for $M_\alpha$ completes the proof.
\end{proof}

The limiting law \eqref{eq:radial-law} makes precise what remains of the
classical ballistic line after the fractional clock is introduced. For
$0<\vartheta<1$, let $e_\vartheta$ be the quantile of order $\vartheta$ of
$\mathcal{E}_1$, and let $r_\vartheta(t)$ be the corresponding radial quantile,
that is,
\[
\mathbb P(R_t\le r_\vartheta(t))=\vartheta.
\]
Since $M_\alpha$ is positive and continuous, $e_\vartheta$ is unique. Moreover, $Z_{\mathbb H}(t,r)>0$ for $r>0$, so the radial distribution of $R_t$ is continuous and strictly increasing; hence $r_\vartheta(t)$ is unique as well.

\begin{corollary}[Radial quantiles]
\label{cor:radial-quantiles}
For every $0<\vartheta<1$,
\begin{equation}
\frac{r_\vartheta(t)}{t^\alpha}
\longrightarrow
(n-1)e_\vartheta,
\qquad t\to\infty.
\label{eq:radial-quantiles}
\end{equation}
In particular, if $0<\vartheta_1<\vartheta_2<1$, then
\[
r_{\vartheta_2}(t)-r_{\vartheta_1}(t)
\sim
(n-1)(e_{\vartheta_2}-e_{\vartheta_1})t^\alpha.
\]
\end{corollary}

\begin{proof}
The limiting distribution in \eqref{eq:radial-law} has a continuous and
strictly increasing distribution function. Convergence of quantiles is then
a standard consequence of weak convergence.
\end{proof}

For the classical heat flow, the mass concentrates in an annulus of radius
$(n-1)s$ and width $O(s^{1/2})$, so the relative width tends to zero.
The preceding corollary shows a different phenomenon for $0<\alpha<1$:
the annulus bounded by two fixed radial quantiles has both radius and width
of order $t^\alpha$. The random inverse clock thus turns the
deterministic ballistic line in operational time into a macroscopically random
radial position. This statement does not mean ballistic propagation in
physical time when $\alpha<1$; the ballistic behavior belongs to the hyperbolic motion in
operational time, and the resulting physical scale is $t^\alpha$.

\subsection{Mean square displacement and the Euclidean comparison}
\label{subsec:properties-msd}

Define the intrinsic mean square displacement by
\begin{equation}
\operatorname{MSD}_\alpha(t)
:=
\mathbb E[R_t^2]
=
\int_{\mathbb H^n}
 d(o,x)^2 Z_{\mathbb H}(t,d(o,x))\,\mathrm d\mu(x).
\label{eq:msd-definition}
\end{equation}
For classical hyperbolic Brownian motion,
\eqref{eq:radial-laplacian} gives
\[
\Delta_{\mathbb H^n}r^2
=
2+2(n-1)r\coth r.
\]
Since $r\coth r\to1$ as $r\downarrow0$, we obtain
\[
m_2(s)\sim2ns,
\qquad s\downarrow0.
\]
Together with \eqref{eq:radial-moment-subordination} and the global estimate
$m_2(s)\le C_n(s+s^2)$, this allows us to pass to the short-time limit.
The large-time behavior is the case $p=2$ of
\eqref{eq:radial-moments}.

\begin{proposition}[Geometric crossover of the MSD]
\label{prop:msd-crossover}
One has
\begin{equation}
\operatorname{MSD}_\alpha(t)
\sim
\begin{cases}
\dfrac{2n}{\Gamma(1+\alpha)}t^\alpha,
& t\downarrow0,\\[3mm]
\dfrac{2(n-1)^2}{\Gamma(1+2\alpha)}t^{2\alpha},
& t\to\infty.
\end{cases}
\label{eq:msd-crossover}
\end{equation}
\end{proposition}

\begin{proof}
By \eqref{eq:radial-moment-subordination}, for $p=2$,
\[
\operatorname{MSD}_\alpha(t)
=
\int_0^\infty m_2(t^\alpha z)M_\alpha(z)\,\mathrm dz.
\]
As $t\downarrow0$, we use $m_2(s)\sim2ns$ and, after division by
$t^\alpha$, the explicit domination
$m_2(t^\alpha z)/t^\alpha\le C_n(z+z^2)$. The moment formula for
$M_\alpha$ gives
\[
\int_0^\infty zM_\alpha(z)\,\mathrm dz
=\frac{1}{\Gamma(1+\alpha)}.
\]
This proves the first line. The second is
\eqref{eq:radial-moments} with $p=2$.
\end{proof}

In $\mathbb R^n$, with the same normalization of the Laplacian, if
$(B_s^{\mathbb R^n})_{s\ge0}$ denotes Euclidean Brownian motion, then
$\mathbb E|B_s^{\mathbb R^n}|^2=2ns$. By subordination,
\[
\operatorname{MSD}^{\mathbb R^n}_\alpha(t)
=
\frac{2n}{\Gamma(1+\alpha)}t^\alpha
\qquad\text{for every }t>0.
\]
Thus the $t^\alpha$ regime is locally Euclidean, whereas the change in the
MSD exponent from $\alpha$ to $2\alpha$ at large times is a genuine consequence
of hyperbolic geometry and its ballistic radial escape. In particular, for $\alpha<1/2$
the intrinsic MSD remains sublinear in physical time, for $\alpha=1/2$ it is
asymptotically linear, and for $\alpha>1/2$ it is superlinear. This last
possibility reflects the geometry and not a loss of the subdiffusive character
of the inverse stable clock.

\subsection*{Data Availability} No datasets were generated or analyzed in this study.

\printbibliography

\end{document}